\documentclass[10pt,a4paper]{article}

\usepackage{amsmath}
\usepackage{amscd}
\usepackage{amssymb}
\usepackage{indentfirst}
\usepackage{latexsym}
\usepackage{multicol}
\usepackage{theorem}

\newcommand{\pf}{{\noindent\it Proof.~}}
\newcommand{\Ov}[1]{\overline{#1}}

\newcommand{\Un}[1]{\underline{#1}}

\newcommand{\vy}{Y_{k}}
\newcommand{\vf}{{\bf F}}
\newcommand{\vQ}{\vc{Q}}
\newcommand{\vC}{{ C}}

\newcommand{\vw}{\omega}

\newcommand{\vr}{\varrho}

\newcommand{\vt}{\vartheta}
\newcommand{\vu}{\vc{u}}
\newcommand{\vd}{\vc{d}_{k}}
\newcommand{\vc}[1]{{\bf #1}}

\newcommand{\Div}{\operatorname{div}}
\newcommand{\Grad}{\nabla}

\newcommand{\pt}{\partial_{t}}
\newcommand{\ptb}[1]{\partial_{t}(#1)}
\newcommand{\Dt}{\frac{{\rm d}} {{\rm dt}}}
\newcommand{\tn}[1]{\mbox {\F #1}}
\newcommand{\dx}{{\rm d} {x}}

\newcommand{\dt}{{\rm d} t }

\newcommand{\dxdt}{\dx \ \dt}
\newcommand{\intO}[1]{\int_{\Omega} #1 \ \dx}

\newcommand{\intOB}[1]{\int_{\Omega} \left( #1 \right) \ \dx}

\newcommand{\intT}[1]{\int_0^T #1 \ \dt}
\newcommand{\intTO}[1]{\int_0^T\!\!\!\! \int_{\Omega} #1 \ \dxdt}

\newcommand{\intTOB}[1]{ \int_0^T\!\!\!\! \int_{\Omega} \left( #1 \right) \ \dxdt}
\newcommand{\sumkN}[1]{\sum_{k=1}^{n} #1}

\newcommand{\D}{{\vc{D}}}
\newcommand{\vS}{{\bf{S}}}

\newcommand{\ep}{\varepsilon}

\font\F=msbm10 scaled 1200
\newcommand{\R}{\mathbb{R}}

\newtheorem{Theorem}{Theorem}

\newtheorem{Lemma}[Theorem]{Lemma}
\newtheorem{Definition}[Theorem]{Definition}
\newtheorem{Remark}[Theorem]{Remark}

\title{Mixtures: sequential stability of variational entropy solutions}
\author{Ewelina Zatorska${}^{1,2}$}
\date{}

\begin{document}
\maketitle
\begin{center}
{\small ${}^1$ {\it CMAP, Ecole Polytechnique, 91128 Palaiseau cedex, France}\\
${}^2$ {\it IAMM, University of Warsaw, ul.Banacha 2, 02-097 Warszawa, Poland}}
\bigskip

\end{center}
\noindent{\bf Abstract:} The purpose of this work is to analyze the mathematical model governing motion of $n$-component, heat conducting reactive mixture of compressible gases.
We prove sequential stability of weak variational entropy solutions when the state equation essentially depends on the species concentration and the species diffusion fluxes depend on gradients of partial pressures. Of crucial importance for our analysis is the fact that viscosity coefficients vanish on vacuum and the source terms enjoy the admissibility condition dictated by the second law of thermodynamics.

\section{Introduction}
To describe the flow of $n$-component chemically reacting compressible mixture, we will use the full Navier-Stokes-Fourier (NSF) system coupled with the set of $n$ reaction-diffusion equations for the species:
\begin{equation}\label{ch4:1.1}
        \begin{array}{lll}
            \vspace{0.2cm}
        \pt\vr+\Div (\vr \vu) = 0& \mbox{in}& (0,T)\times\Omega,\\
            \vspace{0.2cm}
        \ptb{\vr\vu}+\Div (\vr \vu \otimes \vu) - \Div \vS+ \Grad \pi =\vc{0}& \mbox{in}& (0,T)\times\Omega,\\
            \vspace{0.2cm}
        \ptb{\vr E}+\Div (\vr E\vu )+\Div(\pi\vu) +\Div{\vQ}-\Div (\vS\vu)=0\cdot\vu\quad{}& \mbox{in}& (0,T)\times\Omega,\\
            \vspace{0.2cm}
        \pt{\vr_k}+\Div (\vr_{k} \vu)+ \Div (\vf_k)  =  \vr\vt\vw_{k},\quad k\in\{1,...,n\}& \mbox{in}& (0,T)\times\Omega.
            \end{array}
    \end{equation}
These equations express the physical laws of conservation of mass, momentum, total energy and the balances of species mass, respectively. 

Here, $\vu:\R^3\to\R^3$ is the velocity field, $\vr:\R^3\to\R$ denotes the total mass density being  a sum of species densities $\vr_k$, $k\in\{1,\ldots,n\}$. The last unknown quantity is the  temperature $\vt:\R^3\to\R$ which appears implicitly in all the equations of \eqref{ch4:1.1} except for the continuity equation. Next, $\vS$ denotes the viscous stress tensor, the  internal pressure is denoted by $\pi$, $E$ is the total energy per unit mass, $\vQ$ stands for the heat flux, $\vf_k$, $k\in\{1,\ldots,n\}$ denote the species diffusion fluxes and $\vw_k$, $k\in\{1,\ldots,n\}$ are the chemical source terms, also termed the species production rates.

In \eqref{ch4:1.1}, $t$ denotes the time, $t\in(0,T)$ and the length of time interval  $T$ is assumed to be arbitrary large, but finite.
The space domain $\Omega$ is a periodic box $\mathbb{T}^3$. 
The vectors belonging to the physical space $\R^3$ as well as the tensors are printed in the boldface style. 

\noindent The species mass conservation equations can be equivalently written in terms of species mass fractions:
$$ \ptb{\vr Y_k}+\Div (\vr Y_{k} \vu)+ \Div (\vf_k)  =  \vr\vt\vw_{k},$$
where $Y_k$, $k\in\{1,\ldots,n\}$ are defined by 
$Y_k=\frac{\vr_k}{\vr}$ and they satisfy:
\begin{equation}\label{chF:massconserv}
\sumkN Y_k=1.
\end{equation}
We remark that we will freely switch from one notation to the other using the species unknowns $(\vr,\vr_1,\ldots,\vr_n)$ or equivalently $(\vr,Y_1,\ldots, Y_n)$.

The global-in-time existence of solutions for system \eqref{ch4:1.1}, supplemented with physically relevant constitutive relations, was established by Giovangigli [24, chapter 9, Theorem 9.4.1]. He considered, for instance, generic matrix $C_{kl}$ relating the diffusion deriving forces $\vd$ to the species diffusion fluxes
    \begin{equation}\label{eq:diff1}
        \vf_{k}=-\sum_{l=1}^{n} {C}_{k,l}\vc{d}_l~+~ \mbox{Soret\ effect}, \quad k=1,...n,
    \end{equation}
which was singular $C {Y}=0$, $ Y=(Y_1,\ldots,Y_n)^T$,  $Y_k=\frac{\vr_k}{\vr}$, $k=1,\ldots,n$, and  not symmetric in general.
The result holds provided the initial data are sufficiently close to an equilibrium state and our main motivation is to extend it for the case of arbitrary large data. It should be however emphasized that, in comparison with \cite{VG}, many simplifications are adopted in the present model. We concentrate only on diffusion effects due to the mole fraction and pressure gradients and restrict to particular form of the matrix $C$. 
%

Due to mass conservation, the equations for species must sum to continuity equation which is hyperbolic. Moreover, the highly complex structure of diffusion fluxes causes that any standard approach for parabolic systems can not be applied unhindered. Nevertheless, it has been observed that even for strongly coupled systems of PDEs the concept of entropy can still be used in order to derive appropriate compactness results \cite{CJ,GZ}. 
In our case, it is possible provided the matrix $D_{kl}=\frac{C_{kl}}{\vr Y_{k}}$, $k,l=1,\ldots,n$ is symmetric and coercive on the hyperplanes which do not contain the vector ${Y}$. This assumption corresponds to the non-negativity of entropy production rate associated with diffusive process and has been postulated for example by Waldmann \cite{Wal62}. It turns out that such a description captures quite accurately lot of practical applications if only the species deriving forces $\vd$ are well-fitting. 

The most exhaustively studied approximation is the Fick law, widely used especially in the simplified models: two species kinetics, 1D geometry, irreversible or isothermal reactions   (cf. \cite{CHT,DT2,D,EZ,Z}). It states that the diffusion flux of a single species depends only on the gradient of its concentration. Hence, it does not take into account the cross-effects that play an important role in the multicomponent flows. As far as the latter are concerned, the  issue of global-in-time existence of weak variational entropy solutions was investigated by  { Feireisl, Petzeltov\'a and Trivisa} \cite{FTP}. They generalized the proof from \cite{FN} to the case of chemically reacting flows, when there is no interaction among the species diffusion fluxes and the pressure does not depend on the chemical composition of the mixture. The presence of the species concentration in the state equation would result in the entropy production rate which may fail to be non-negative. This in turn interferes with obtaining  the fundamental a-priori estimates being the corner stone of the analysis presented in \cite{FTP}.  A similar problem was apparent in the approximation considered by Frehse, Goj and M\'alek \cite{FGM05a, FGM05b} who proved existence and uniqueness of solutions to the steady Stokes-like system for a mixture of two (non-reactive) fluids. In their case, neglecting some nonlinear interaction terms in the source of momentum caused that the basic energy equality was no longer preserved. 

Significantly  less is known for models with general diffusion. The local-in-time well posedness of the Maxwell-Stefan equations with multicomponent diffusion  in the isobaric, isothermal case is due to Bothe \cite{B2010}. More recently, Mucha, Pokorn\'y and Zatorska \cite{MPZ} proved that the system of $n$ diffusion-reaction equations with particular form of diffusion matrix $C$ admits a global-in-time weak solution, provided the additional regularity of the total density is available. 
This extra information can be deduced from the new kind of estimate found for the  Korteweg system by Bresch, Desjardins and Lin \cite{BDL}.  It also holds for the compressible Navier-Stokes equations \cite{BDG} if the viscosity coefficients $\mu$ and $\nu$ are density-dependent and the stisfy
$$\nu(\vr)=2\vr \mu'(\vr)-2\mu(\vr).$$
On the other hand, such an assumption leads automatically  to a problem with defining the velocity field on the vacuum zones. It was solved on the level of weak sequential stability of solutions by Mellet and Vasseur  \cite{MV07} who combined the previous results with logarithmic estimate for the velocity. Nevertheless,  construction of sufficiently smooth approximate solutions in this framework remains an open problem. A way to overcome this obstacle in the case of heat-conducting fluids was proposed in a recent paper by Bresch and Desjardins \cite{BD}.  They introduced a modification of the pressure close to the zero Kelvin isothermal curve by which it was possible to show that the sets of presence of vacuum have a zero Lebesgue measure. The same idea was then employed in \cite{EZ2,MPZ1} to prove the global-in-time existence of weak solutions for the isothermally reacting mixture of two gaseous components.
%
%
%
In the present paper we indicate a possible way of generalization of these results to the case of $n$-component heat-conducting mixture.

%
Let us emphasize, that unlike to \cite{BD} we use the framework of weak variational solutions, that is to say, the total energy balance is replaced by the entropy inequality and the global total energy balance. These solutions were introduced by Feireisl \cite{EF1} to study solvability of the NSF system, for which the weak solutions may dissipate more kinetic energy than, if there are any, the classical ones. However, it can be verified that this "missing energy" is equal to $0$ provided the weak solution is regular enough and satisfies the total global energy balance. Moreover, in \cite{FN2012} the authors proved  the weak-strong uniqueness of such solution, meaning that it coincides with the strong solution, emanating from the same initial data, as long as the latter exists. 

The article is organized as follows. In the next section we specify the structural properties for the transport coefficients and postulate several simplifications. In Section \ref{ch4:sect3}, we define the notion of weak variational solutions and state the main result of the paper. The key a priori estimates are derived in Section \ref{ch4:sect4} together with some further estimates and positivity of the absolute temperature. Finally, the last step of the proof of Theorem \ref{ch4:Theorem:1} - the limit passage - is performed in Section \ref{ch4:sect6}.
%
%

\section{Constitutive relations and main hypothesis}
Let us now  supplement system \eqref{ch4:1.1} with a set of expressions determining the form of thermodynamical functions and transport fluxes in terms of  macroscopic variable gradients in the spirit of \cite{VG}, Chapter 2. \\

\noindent{\bf Thermal equation of state.} We consider the pressure $\pi=\pi(\vr,\vt,Y)$, which can be decomposed into
\begin{equation}\label{chF:decompp}
\pi=\pi_m+\pi_c,
\end{equation}
where the latter component depends solely on the density and it corresponds to the barotropic process of viscous gas. It is the only non-vanishing component of the pressure when temperature goes to zero, thus will be termed a "cold" pressure. We assume that $\pi_c$ is a continuous function satisfying the following growth conditions
\begin{equation}\label{ch4:coldp}
\pi'_{c}(\vr)=\left\{
\begin{array}{cl}
c_1\vr^{-\gamma^- -1}&\mbox{for} \ \vr\leq 1,\\
c_2\vr^{\gamma^+-1}&\mbox{for} \ \vr>1,
\end{array}\right.
\end{equation}
for positive constants $c_1,c_2$ and $\gamma^-,\gamma^+>1$.
It was proposed in \cite{BD} to encompass plasticity and elasticity effects of solid materials, for which low densities may lead to negative pressures. By this modification the compactness of velocity at the last level of approximation can be obtained without requiring more a priori regularity than expected from the usual energy approach  \cite{MV07}.\\

\noindent The first component $\pi_{m}=\pi_m(\vr,\vt,Y)$ is the classical molecular pressure of the mixture which is determined through the {\it Boyle law} as a sum of partial pressures $p_k$:
$$\pi_{m}(\vr,\vt,Y)=\sum_{k=1}^{n}p_{k}(\vt,\vr_{k})=\sum_{k=1}^{n}\frac{\vt\vr_{k}}{m_{k}},$$
with $m_{k}$ the molar mass of the $k$-th species.\\
Likewise the pressure, the internal energy $e=e(\vr,\vt,Y)$ can be decomposed into
\begin{equation}\label{chF:enerdef}
e=e^{st}+e_{m}+e_{c},\quad e^{st}(Y)=\sumkN Y_{k}e^{st}_{k},\quad e_m(\vt)= c_v\vt,
\end{equation}
where $e^{st}_{k}={\it const.}$ is the formation energy of the $k$-th species, while $c_{v}$ is the constant-volume specific heat, which is assumed to be the same for each of the species. The "cold" components of the internal energy $e_c=e_c(\vr)$ and pressure $\pi_c$ are related through the following equation of state:
\begin{equation}\label{chF:intener}
\vr^2\frac{\mbox{d}e_c(\vr)}{\mbox{d}\vr}=\pi_c(\vr).
\end{equation}
The last relation is a consequence of the second law of thermodynamics which postulates the existence of a state function called the entropy.
\vspace{0.5cm}

\noindent{\bf The entropy equation.}
The entropy of a thermodynamical system is defined (up to an additive constant) by the differentials of energy, total density, and species mass fractions via the Gibbs relation:
\begin{equation}\label{chF:entropdiff}
\vt \mathrm{D} s=\mathrm{D} e+\pi\mathrm{D}\left({1 \over \vr}\right)-\sumkN g_{k}\mathrm{D} \vy,
\end{equation}
where $\mathrm{D}$ denotes the total derivative with respect to the state variables $\{\vr,\vt,Y\}$; whereas $g_k$ are the Gibbs functions
\begin{equation}\label{chF:defg}
g_{k}=h_{k}-\vt s_{k}.
\end{equation}
Here, $h_k=h_k(\vt)$ denotes the specific enthalpy and $s_k=s_k(\vt,\vr_{k})$ is the specific entropy of the $k$-th species
$$ h_k(\vt)=e^{st}_{k}+c_{pk}\vt,$$
\begin{equation}\label{chF:entropycv}
s_k(\vt,\vr_k)=s_{k}^{st}+c_{v}\log{{\vt}}+{1\over m_{k}}\log{ {m_{k}}\over{\vr_k}},
\end{equation}
where $s^{st}_{k}=const.$ denotes the formation entropy of $k$-th species.\\
The constant-volume ($c_{vk}$) and constant-pressure ($c_{pk}$) specific heats for the $k$-th species  are constants related by 
\begin{equation}\label{chF:cpcv}
c_{pk}=c_{v}+\frac{1}{m_k}.
\end{equation}
In accordance with \eqref{chF:enerdef},  \eqref{chF:intener} and \eqref{chF:entropdiff}, the specific entropy of the mixture can be expressed as a weighted sum of the species specific entropies
\begin{equation}\label{chF:cotos}
s=\sumkN\vy s_{k}
\end{equation}
and is governed by the following equation
\begin{equation}\label{chF:entropy}
\ptb{\vr s}+\Div(\vr s\vu)+\Div\left( \frac{\vQ}{\vt}-\sumkN \frac{g_{k}}{\vt}\vf_{k}\right)=\sigma,
\end{equation}
where $\sigma$ is the entropy production rate
\begin{equation} \label{chF:sigma}
\sigma={\vS:\Grad\vu\over\vt}-{\vQ\cdot\Grad\vt\over{\vt^{2}}}-\sumkN\vf_{k}\cdot\Grad\left({g_{k}\over \vt}\right)-{{\sumkN g_{k}\vr\vt\vw_{k}}\over\vt}.
\end{equation}
By virtue of the second law of thermodynamics, the entropy production rate must be non-negative for any admissible process.

\vspace{0.5cm}

\noindent {\bf The stress tensor.} The viscous part of the stress tensor obeys the {\it Newton rheological law}, namely:
\begin{equation}\label{chF:Stokes}
\vS= 2\mu\D(\vu)+\nu\Div \vu\tn{I},
\end{equation}
where $\D(\vu)=\frac{1}{2}\left(\Grad \vu+(\Grad \vu)^{T}\right)$ and $\mu=\mu(\vr),\nu=\nu(\vr)$
are $C^{1}(0,\infty)$ functions related by
\begin{equation}\label{ch4:9}
\nu(\vr)=2\vr \mu'(\vr)-2\mu(\vr),
\end{equation}
which is strictly a mathematical constraint allowing to obtain better regularity of $\vr$, \cite{BD2}.\\
In addition, they enjoy the following bounds
\begin{equation}\label{ch4:10}
\begin{array}{c}
\Un{\mu'}\leq\mu'(\vr)\leq{1 \over \Un{\mu'}},\quad \mu(0)\geq 0,\\
|\nu'(\vr)|\leq{1\over \Un{\mu'}}\mu'(\vr),\\
\Un{\mu'} \mu(\vr)\leq 2\mu(\vr)+3\nu(\vr)\leq{1\over \Un{\mu'}}\mu(\vr),
\end{array}
\end{equation}
for some positive constant  $\Un{\mu'}$,
\begin{Remark}
The assumption $\Un{\mu'}\leq\mu'(\vr)$ is not optimal, but it makes the proof much easier, however the bound from above is essential in order to get integrability of several important quantities. For further discussion on this topic we refer to \cite{MV07}.
\end{Remark}


\vspace{0.5cm}

\noindent {\bf The species diffusion fluxes.} Following \cite{VG}, Chapter 2, Section 2.5.1, we postulate that the diffusion flux of the $k$-th species is given by
    \begin{equation}\label{chF:eq:diff}
        \vf_{k}=-C_{0}\sum_{l=1}^{n} {C}_{kl}\vc{d}_l, \quad k=1,...n,
    \end{equation}
where $C_{0},\ C_{kl}$ are multicomponent flux diffusion coefficients and $\vd$ is the species $k$ diffusion force specified, in the absence of external forces, by the following relation
        \begin{equation}\label{chF:eq:}
        \vd=\Grad\left({p_{k}\over \pi_m}\right)+\left({p_{k}\over \pi_m}-{\vr_{k}\over \vr}\right)\Grad\log{\pi_m}.
    \end{equation}
We restrict to an exact form of the flux diffusion matrix $\vC$ compatible with the set of mathematical assumptions postulated in \cite{VG}, Chapters 4, 7 (see also references therein). The prototype example is the same as studied in \cite{MPZ}, namely 
    \begin{equation}\label{ch4:prop}
        \vC =\left(
            \begin{array}{cccc}
                Z_{1} & -Y_{1} & \ldots & -Y_{1}\\
                -Y_{2} & Z_{2}  & \ldots & -Y_{2}\\
                \vdots & \vdots & \ddots & \vdots \\
                -Y_{n} & -Y_{n}& \ldots & Z_{n}
            \end{array}
        \right),
    \end{equation}
where  $Z_{k}=\sum_{{i=1} \atop {i\neq k}}^{n} Y_{i}$.\\
Concerning the diffusion coefficient $C_0$, we assume that it is a continuously differentiable function of $\vt$ and $\vr$  and that there exist positive constants $\Un{\vC_0},\ \Ov{\vC_0} $ such that
 \begin{equation}\label{ch4:assC}
\Un{C_0}\vr(1+\vt)\leq\vC_{0}\leq \Ov{C_0}\vr(1+\vt).
\end{equation}
\begin{Remark}
One of the main consequences of \eqref{ch4:prop} is that
\begin{equation} \label{ch4:zerototdiff}
    \sum_{k=1}^{n}\vf_{k}=0.
\end{equation}
\end{Remark}
\begin{Remark}\label{ch4:remF}
Note that \eqref{ch4:prop} also implies  that the vector of species diffusion forces ${\vc{d}}=(\vc{d}_1,\ldots,\vc{d}_n)^T$ is an eigenvector of the matrix $\vC$ corresponding to the eigenvalue $1$ 
     and we recast that
    \begin{equation}\label{ch4:difp}
        \vf_{k}=-\vC_0\sum_{l=1}^nC_{kl}\vc{d}_l=-\vC_0\vc{d}_k=-\frac{\vC_0}{\pi_m}\left(\Grad p_{k}-Y_{k}\Grad \pi_m\right)=-\frac{\vC_0}{\pi_m}\sum_{l=1}^{n}\vC_{kl}\Grad p_{l}.
    \end{equation}
 \end{Remark}
 \vspace{0.5cm}

\vspace{0.2cm}

\noindent {\bf The heat flux.} It is a sum of two components
    \begin{equation}\label{chF:cotoQ}
        \vQ=\sumkN h_k \vf_{k}-\kappa\Grad\vt,
    \end{equation}
where the first term describes transfer of energy due to the species molecular diffusion, whereas the second term represents the {\it Fourier law} with the thermal conductivity coefficient $\kappa=\kappa(\vr,\vt)$. It is assumed to be a $C^1([0,\infty)\times [0,\infty))$ function which satisfies
\begin{equation}\label{ch4:cotok}
\Un{\kappa_0}(1+\vr)(1+\vt^\alpha)\leq\kappa(\vr,\vt)\leq \Ov{\kappa_0}(1+\vr)(1+\vt^\alpha).
\end{equation}
In the above formulas $\Un{\kappa_0}, \Ov{\kappa_0},\alpha$ are positive constants and $\alpha\geq2$. \\


\noindent{\bf{The species production rates}}. 
We assume that $\vw_k$ are  continuous functions of $Y$ bounded from below and above by the positive constants $\Un{\vw}$ and $\Ov{\vw}$ 
    \begin{equation}\label{ch4:wform}
        -\Un{\vw}\leq\vw_{k}(Y)\leq\Ov{\vw},\quad  \quad \mbox{ for\ all}\ k=1,\ldots,n;
    \end{equation}
we also suppose that
    \begin{equation}\label{ch4:wform1}
    \sum_{k=1}^{n}\vw_{k}=0,\quad\mbox{and}\quad
    \vw_{k}(Y)\geq 0\quad\mbox{ whenever}\ \ Y_{k}=0.
    \end{equation}
Another restriction is dictated by the second law of thermodynamics, which asserts in particular that $\vw_{k}$ must enjoy the following  condition
\begin{equation}\label{ch4:admiss0}
\sumkN g_{k}\vr\vw_{k}\leq0.
\end{equation}



\noindent{\bf Initial data.} The choice of quantities describing the initial state of system \eqref{ch4:1.1} is dictated by the weak formulation of the problem, which is  specified in Definition \ref{ch4:def1} below. We take
\begin{equation}\label{ch4:initial}
\begin{array}{c}
\vspace{0.2cm}
\vr(0,\cdot)=\vr^{0},\ \vr\vu(0,\cdot)=(\vr\vu)^{0},\ \vr s(0,\cdot)=(\vr s)^{0},\ \intO{\vr E(0,\cdot)}=\intO{(\vr E)^0},\\
\vr_k(0,\cdot)=\vr_k^{0}, \ \mbox{ for}\ k=1,\ldots,n,\ \mbox{in}\ \Omega.
\end{array}
\end{equation}
In addition, we assume that $\vr^0$ is a nonnegative measurable function such that
\begin{equation}\label{ch4:initial_mass}
\intO{\vr^0}=M^0,\quad \intO{\frac{1}{\vr^{0}}\left|\Grad\mu\left(\vr^{0}\right)\right|^{2}}<\infty,
\end{equation}
and the initial densities of species satisfy
\begin{equation}\label{ch4:initialab}
0\leq\vr_k^0(x),\  k=1,\ldots,n,\quad \sumkN\vr_k^0(x)=\vr^0(x), \ \mbox{a.e.\ in\ }\Omega.
\end{equation}
Further, we call the initial temperature $\vt^0$ a measurable function such that
$$
\vt^0(x)>0\  \mbox{a.e.\ in\ }\Omega,\quad\vt^0\in W^{1,\infty}(\Omega)$$
and the following compatibility condition is satisfied
\begin{equation}\label{ch4:comp_teta}
(\vr s)^0=\vr^0s(\vt^0,\vr_1^0,\ldots,\vr_n^0),\quad (\vr s)^0\in L^1(\Omega).
\end{equation}
Finally, we require that the initial distribution of the momentum is such that
$$(\vr\vu)^0=0\ \mbox{a.e.\ on\ }\{x\in\Omega:\vr^0(x)=0\}\quad\mbox{and}\quad \intO{\frac{\left|(\vr\vu)^0\right|^2}{\vr^0}}<\infty$$
and the global total energy at the initial time is bounded
\begin{equation}\label{ch4:glob_int_ener}
\intO{(\vr E)^0}=\intOB{\frac{\left|(\vr\vu)^0\right|^2}{2\vr^0}+\vr^0 e(\vr^0,\vt^0,\vr_1^0,\ldots,\vr_n^0)}<\infty.
\end{equation}

\section{Weak formulation and main result}\label{ch4:sect3}
\noindent In what follows we define a notion of weak variational entropy solutions to system \eqref{ch4:1.1} and then we formulate our main result.

\begin{Definition}\label{ch4:def1}
We will say  $\{\vr,\vu,\vt,\vr_1,\ldots,\vr_n\}$ is a weak variational entropy solution provided the following integral identities hold.\\
1. The continuity equation
 \begin{equation}\label{ch4:weakcont}
\intO{\vr_0\phi(0,x)}+\intTO{\vr\pt\phi}+\intTO{\vr\vu\cdot\Grad\phi}=0
\end{equation}
is satisfied for any smooth function $\phi(t,x)$ such that $\phi(T,\cdot)=0.$\\
2. The balance of momentum
\begin{multline}\label{ch4:weakmom}
\intO{(\vr\vu)^0\cdot\phi(0,x)}+\intTOB{{\vr}\vu\cdot\pt\phi+{\vr}\vu\otimes\vu:\Grad\phi}\\
+\intTO{\pi(\vr,\vt,Y)\Div\phi}-\intT{\vS:\Grad\phi}=0
\end{multline}
holds for any smooth test function
$\phi(t,x)$ such that $\phi(T,\cdot)=0$. \\
3. The entropy equation
\begin{multline}\label{ch4:weakentr}
\intO{\vr^0 s(\vt^0,\vr_1^0,\ldots,\vr_n^0)\phi(0,x)}+\intTO{\vr s\pt\phi}+\intTO{\vr s\vu\cdot\Grad\phi}\\
+\intTO{\left(\frac{\vQ}{\vt}-\sumkN\frac{g_k}{\vt}\vf_k\right)\cdot\Grad\phi}
+{\left\langle\sigma,\phi \right\rangle}=0
\end{multline}
is satisfied for any smooth function $\phi(t,x)$, such that $\phi\geq0$ and $\phi(T,\cdot)=0$,
where $\sigma\in{\cal M}^+([0,T]\times\Omega)$ is a nonnegative measure such that
$$\sigma\geq
{\frac{\vS:\Grad\vu}{\vt}
-{\vQ\cdot\Grad\vt\over{\vt^{2}}}-\sumkN\frac{\vf_k}{m_k}\cdot\Grad\left(\frac{g_k}{\vt}\right)
-{\sumkN g_k \vr \vw_k}}.$$
4. The global balance of total energy
\begin{equation}\label{ch4:weaktoten}
\intO{(\vr E)^0} ~\phi(0)+\intT{\intO{\vr E}~\pt\phi(t)}=0
\end{equation}
holds for any smooth function $\phi(t)$, such that $\phi(T)=0$.\\
5. The weak formulation of the mass balance equation for the $k$-th species
\begin{multline}\label{ch4:weakspec}
\intO{\vr_k^0\cdot\phi(0,x)}+\intTOB{\vr_k\pt\phi+\vr_k\vu\cdot\Grad\phi}
-\intTO{\vf_k\cdot\Grad\phi}\\
=\intTO{\vr\vt\vw_k\phi}
\end{multline}
$k=1,\ldots,n,$ is satisfied  for any smooth test function
$\phi(t,x)$ such that $\phi(T,\cdot)=0$

In addition we require that 
\begin{equation}\label{ch4:weaksum}
\vr,\vr_k\geq 0,\ k=1,\ldots,n,\quad \sumkN\vr_k=\vr, \quad \mbox{and}\quad\vt>0,\quad\mbox{a.e.\ on\ }(0,T)\times\Omega.
\end{equation}
\end{Definition}

In this definition the usual weak formulation of the energy equation is replaced by the weak formulation of the entropy inequality and the global total energy balance \eqref{ch4:weakentr}+\eqref{ch4:weaktoten}. Note however that the entropy production rate has now only a meaning of non-negative measure which is bounded  from below by the classical value of $\sigma$. Nevertheless, a simple calculation employing the Gibbs formula \eqref{chF:entropdiff} shows that whenever the solution specified above is sufficiently regular, both formulations are equivalent. In particular the entropy inequality \eqref{ch4:weakentr} changes into equality, see f.i. Section 2.5 in \cite{F07}.

We are now in a position to formulate the main result of this paper.
\begin{Theorem}\label{ch4:Theorem:1}
 Assume that the structural hypotheses (\ref{chF:decompp}-\ref{ch4:admiss0})  are satisfied.\\ 
Suppose that $\{\vr_N,\vu_N,\vt_N,\vr_{k,N}\}_{N=1}^{\infty}$, $k=1,\ldots,n$ is a sequence of smooth solutions to (\ref{ch4:1.1})  with the initial data
\begin{equation*}
\begin{array}{c}
\vspace{0.2cm}
\vr_N(0,\cdot)=\vr_N^{0},\ (\vr\vu)_N(0,\cdot)=(\vr\vu)_N^0,\ (\vr s)_N(0,\cdot)=(\vr s)_N^0, \ \intO{(\vr E)_N(0,\cdot)}=\intO{(\vr E)_N^0},\\
\vr_{k,N}(0,\cdot)=\vr_{k,N}^{0}, \ \mbox{ for}\ k=1,\ldots,n,\ \mbox{in}\ \Omega.
\end{array}
\end{equation*}
satisfying (\ref{ch4:initial_mass}-\ref{ch4:glob_int_ener}), moreover
\begin{equation}\label{ch4:initial2}
\begin{array}{c}
\vspace{0.2cm}
\inf_{x\in\Omega}\vr_N^{0}(x)>0,\quad \inf_{x\in\Omega}\vt_N^{0}(x)>0,\\
\vspace{0.2cm}
\vr_N^{0}\rightarrow\vr^{0},\ (\vr\vu)_N^0\to(\vr\vu)^0,\ (\vr s)_N^0\to(\vr s)^0,\ {(\vr E)_N^0}\to{(\vr E)^0},\ \vr_{k,N}^{0}\to\vr_{k}^{0}\quad in\ L^{1}(\Omega).
\end{array}
\end{equation}
Then, up to a subsequence, $\{\vr_N,\vu_N,\vt_N,\vr_{k,N}\}$ converges to the weak solution of problem (\ref{ch4:1.1}) in the sense of Definition \ref{ch4:def1}. 
\end{Theorem}
The proof of this theorem can be divided into two main steps. The first one is dedicated to derivation of the energy-entropy estimates which are obtained under assumption that all the quantities are sufficiently smooth. The next step is the limit passage, it consists of various integrability lemmas  combined with condensed compactness arguments.

\section{A priori estimates}\label{ch4:sect4}
In this section we present the a priori estimates for $\{\vr_N,\vu_N,\vt_N,\vr_{k,N}\}_{N=1}^\infty$ which is a sequence of smooth functions solving \eqref{ch4:1.1}. As mentioned above, assuming smoothness of solutions, we expect that all the natural features of the system can be recovered. The following estimates are valid for each $N=1,2,\ldots$ but we skip the subindex when it does not lead to any confusion.
\subsection{Estimates based on the maximum principle}
To begin, observe that the total mass of the fluid is a constant of motion, meaning
\begin{equation}\label{ch4:M0}
\intO{\vr(t,x)}=\intO{\vr^0}=M_0\quad\mbox{for}\ t\in[0,T].
\end{equation}
Moreover if solution is sufficiently smooth, a classical maximum principle can be applied to the continuity equation in order to show that $\vr_N(t,x)>c(N)\geq0$, more precisely
\begin{equation}\label{ch2:ron0}
\vr_N(\tau,x)\geq \inf_{x\in\Omega}\vr^0_N(x)\exp\left(-\int_{0}^\tau\|\Div\vu_N\|_{L^\infty(\Omega)}{\rm d}t\right),
\end{equation}
in particular $\vr>0$.\\

Next, by a very similar reasoning we can prove non-negativity of $\vt$ on $[0,T]\times\Omega$. 
\begin{Lemma}
Assume that $\vt=\vt_N$ is a smooth solution of \eqref{ch4:1.1}, then
\begin{equation}\label{ch4:pos_temp}
\vt(t,x)> c(N)\geq0 \quad\mbox{for}\ (t,x)\in [0,T]\times\Omega.
\end{equation}
\end{Lemma}
\pf Any solution to \eqref{ch4:1.1} which is sufficiently smooth is automatically a classical solution of the system
\begin{equation}\label{ch4:1.11}
        \begin{array}{l}
            \vspace{0.2cm}
        \pt\vr+\Div (\vr \vu) = 0,\\
            \vspace{0.2cm}
        \ptb{\vr\vu}+\Div (\vr \vu \otimes \vu) - \Div \vS+ \Grad \pi =\vc{0},\\
            \vspace{0.2cm}
       \pt(\vr e)+\Div(\vr e\vu)+\Div\vQ=-\pi\Div\vu+\vS:\Grad\vu,\\
            \vspace{0.2cm}
        \pt{\vr_k}+\Div (\vr_{k} \vu)+ \Div (\vf_k)  =  \vr\vt\vw_{k},\quad k\in\{1,...,n\},\\
            \end{array}
    \end{equation}
 where we replaced the total energy balance by the internal energy balance.
Let us hence reformulate it in order to obtain the equation describing the temperature. Subtracting from the third equation of \eqref{ch4:1.11} 
the component corresponding to the formation energy we obtain
\begin{multline*}
\pt(\vr (e_c+e_m))+\Div(\vr (e_c+e_m)\vu)+\sumkN\Div(c_{pk}\vt\vf_{k})-\Div(\kappa(\vr,\vt)\Grad\vt)\\
=-\pi\Div\vu+\vc{S}:\Grad\vu-\vr\vt\sumkN e^{st}_{k}\omega_{k},
\end{multline*}
where we used the species mass balance equations. Next renormalizing the continuity equation and employing the relation between  $e_{c}$ and $\pi_{c}$ \eqref{chF:intener} we get the temperature equation 
\begin{multline}\label{ch4:tempeq}
\pt(\vr e_{m})+\Div(\vr e_{m}\vu)+\sumkN\Div(c_{pk}\vt\vf_{k})-\Div(\kappa(\vr,\vt)\Grad\vt)\\
=-\pi_{m}\Div\vu+\vc{S}:\Grad\vu-\vr\vt\sumkN e^{st}_{k}\omega_{k},
\end{multline}
where, in accordance with hypotheses \eqref{ch4:10} the second term on the r.h.s. is nonnegative. Consequently, \eqref{ch4:pos_temp} is obtained by application of the maximum principle to the above equation, recalling that $\rm{ess}\inf_{x\in\Omega}\vt_N^{0}(x)>0$. $\Box$\\

The analogous result for partial masses is stated in the following lemma.
\begin{Lemma}\label{ch4:Lemma1}
For any smooth solution of (\ref{ch4:1.1}) we have
    \begin{equation}\label{ch4:dodY}
        \vr_{k}(t,x)\geq0\quad {\rm for}\ (t,x)\in [0,T]\times\Omega,\ k\in\{1,\ldots,n\}.
    \end{equation}
Moreover 
\begin{equation}\label{ch4:sumab}
\sumkN\vr_k(t,x)=\vr(t,x)\quad {\rm for}\  (t,x)\in[0,T]\times\Omega.
\end{equation}

\end{Lemma}
\pf We integrate each of equations of system \eqref{ch4:1.1} over the set $\{\vr_{k}<0\}$. Assuming that the boundary i.e. $\{\vr_{k}=0\}$ is a regular submanifold we obtain
    $$\Dt\int_{\{\vr_{k}<0\}}\vr_{k}~\dx-\int_{\{\vr_{k}=0\}}\frac{\partial p_{k}}{\partial n}~\mbox{d}S_{x}+\int_{\{\vr_{k}=0\}}\frac{\vr_{k}}{\vr}\frac{\partial \pi_m}{\partial n}~    \mbox{d}S_{x}   =\int_{\{\vr_{k}<0\}}\vr\vt\vw_{k}~\dx.$$
Since $\frac{\partial p_{k}}{\partial n}\big|_{\{\vr_{k}=0\}}\geq 0$ and $\vw_{k}\big|_{\{\vr_{k}<0\}}\geq0$ we find
    $$\int_{\{\vr_{k}<0\}}\vr_{k}(T)~\dx\geq\int_{\{\vr_{k}<0\}}\vr_{k}^{0}~\dx=0,$$
thus $|\{\vr_{k}<0\}|=0,$ for every $k=1,\ldots,n.$ When $\{\vr_{k}=0\}$ is not a regular submanifold we construct a sequence $\{\ep_l\}_{l=1}^\infty$ such that $\ep_l\to 0^{+}$ and $\{\vr_{k}=\ep_l\}$ is a regular submanifold and pass with $\ep_l$ to zero.\\
The proof of \eqref{ch4:sumab} follows by subtracting the sum of species mass balances equations from the continuity equation. The smooth solution of the resulting system must be, due to the initial conditions \eqref{ch4:initialab}, equal to $0$ on $[0,T]\times\Omega$. $\Box$\\

As a corollary from this lemma we recover relation \eqref{chF:massconserv}, moreover we have the following estimate
   \begin{equation}\label{ch4:maxY}
        \|Y_{k}\|_{L^{\infty}((0,T)\times\Omega)}\leq 1, \quad k=1\ldots,n.
    \end{equation}

\subsection{The energy-entropy estimates}
The purpose of this subsection is to derive a priori estimates resulting from the energy and entropy balance equations. The difference comparing to estimates obtained in the previous subsection is that now we look for bounds which are uniform with respect to $N$. We start with the following Lemma.

\begin{Lemma}\label{ch4:Lemma3}
Every smooth solution of \eqref{ch4:1.1} satisfies
\begin{equation}
\Dt\intO{\vr\left(\frac{1}{2}|\vu|^{2}+e\right)}=0.
\end{equation}
\end{Lemma}
\pf Integrate the third equation of \eqref{ch4:1.1}  with respect to the space variable and employ the periodic boundary conditions. $\Box$

Assuming integrability of the initial conditions \eqref{ch4:glob_int_ener} the assertion of the above lemma entails several a priori estimates:
\begin{equation}\label{ch4:vrvt}
\begin{array}{c}
\vspace{0.2cm}
\|\sqrt{\vr}\vu\|_{L^{\infty}(0,T;L^{2}(\Omega))}\leq c,\\
\|\vr e_c(\vr)\|_{L^{\infty}(0,T;L^{1}(\Omega))}+\|\vr\vt\|_{L^{\infty}(0,T;L^{1}(\Omega))}+\|\vr\|_{L^{\infty}(0,T;L^{1}(\Omega))}\leq c.\
\end{array}
\end{equation}
It is well known that these natural bounds are not sufficient to prove the weak sequential stability of solutions, not even for the barotropic flow. However, taking into account the form of viscosity coefficients \eqref{ch4:9}, \eqref{ch4:10}, further estimates can be delivered. 
\begin{Lemma}\label{ch4:Lemma4}
For any smooth solution of \eqref{ch4:1.1} we have
\begin{equation}\label{ch4:first}
\Dt \intO{{1\over2}\vr|\vu|^2}+\intO{\vS:\Grad\vu}=\intO{\pi(\vr,\vt,Y)\Div\vu}.
\end{equation}
\end{Lemma}
\pf Multiply the momentum equation by $\vu$ and integrate over $\Omega$. $\Box$\\
The above lemma can not be used to deduce the uniform bounds for the symmetric part of the gradient of $\vu$ immediately. The reason for that is lack of sufficient information for $\vt$, so far we only know \eqref{ch4:vrvt}. 
However, taking into account the form of viscosity coefficients \eqref{ch4:9}, further estimates can be delivered. The following lemma is a modification of the result proved by Bresch and Desjardins for heat conducting fluids \cite{BD}.
%
\begin{Lemma}\label{ch4:Lemma5}
Any smooth solution of \eqref{ch4:1.1} satisfies the following identity
\begin{multline}\label{ch4:both}
\Dt \intO{{1\over2}\vr|\vu+\Grad\phi(\vr)|^2}+{1\over2}\intO{\mu(\vr)|\Grad\vu-\Grad^{T}\vu|^2}=\\
-\intO{\Grad\phi(\vr)\cdot\Grad \pi(\vt,\vr,Y)}+\intO{\pi(\vr,\vt,Y)\Div\vu},
\end{multline}
for $\phi$ such that
$\Grad\phi(\vr)=2\frac{\mu'(\vr)\Grad\vr}{\vr}$.
\end{Lemma}
\pf The rough idea of the proof is the following.
The terms from the l.h.s. of this equality can be evaluated by multiplication of the momentum equation by $\Grad\phi(\vr)$ and the continuity equation by $|\Grad\phi(\vr)|^2$. Then one has to combine these equivalences with the balance of kinetic energy \eqref{ch4:first}
and include \eqref{ch4:9} to see that some unpleasant terms cancel. For more details we refer to \cite{EZ2} or to the original work of Bresch and Desjardins \cite{BD} . $\Box$\\

To control the r.h.s. of \eqref{ch4:first} and \eqref{ch4:both} one needs i.a. to estimate the gradient of $\vt$. To this purpose we take advantage of the entropy balance \eqref{chF:entropy}, we have the following inequality
\begin{Lemma}\label{ch4:Lemma6}
For any smooth solution of \eqref{ch4:1.1} we have
\begin{multline}\label{ch4:estentr}
\intTO{{\vS:\Grad\vu\over\vt}}+\intTO{{\kappa|\Grad\vt|^2\over{\vt^{2}}}}+\intTO{\sumkN\frac{\pi_m\vf_{k}^{2}}{\vC_0\vt\vr_{k}}}\\
-\intTO{{{\sumkN g_{k}\vr\vw_{k}}}}\leq c.
\end{multline}
\end{Lemma}
\pf Combining the third equation of  \eqref{ch4:1.1} with the Gibbs relation \eqref{chF:entropdiff} we derive the entropy equation 
\begin{multline}\label{ch4:ent}
\ptb{\vr s}+\Div(\vr s\vu)+\Div\left( \frac{\vQ}{\vt}-\sumkN \frac{g_{k}}{\vt}\vf_{k}\right)\\
={\vS:\Grad\vu\over\vt}-{\vQ\cdot\Grad\vt\over{\vt^{2}}}-\sumkN\vf_{k}\cdot\Grad\left({g_{k}\over \vt}\right)-{{\sumkN g_{k}\vr\vw_{k}}}.
\end{multline}
Integrating it over space and time we obtain
\begin{equation*}
\intTOB{\frac{\vS:\Grad\vu}{\vt}
-{\vQ\cdot\Grad\vt\over{\vt^{2}}}-\sumkN\vf_k\cdot\Grad\left(\frac{g_k}{\vt}\right)
-{\sumkN g_k \vr \vw_k}}
=\intO{\vr s(T)}-\intO{(\vr s)^{0}},
\end{equation*}
where the l.h.s. can be transformed using \eqref{chF:defg} and \eqref{chF:cotoQ} into 
\begin{multline}\label{ch4:sigma2}
\intTO{{\vS:\Grad\vu\over\vt}}+\intTO{{\kappa|\Grad\vt|^2\over{\vt^{2}}}}-\intTO{\sum_{k=1}^{n}\frac{\vf_{k}} {m_k}\cdot\Grad\log{p_{k}}}\\
-\intTO{{{\sumkN g_{k}\vr\vw_{k}}}}=\intO{\vr s(T)}-\intO{(\vr s)^{0}}.
\end{multline}
The first two terms on the l.h.s. of (\ref{ch4:sigma2}) have a good sign, the same holds for the last one due to \eqref{ch4:admiss0}. Non-negativity of the third one follows from \eqref{ch4:zerototdiff} and  \eqref{ch4:difp}
\begin{equation*}
-\sum_{k=1}^{n}\frac{\vf_{k}} {m_k}\cdot\Grad\left(\log{p_{k}}\right)=-\sum_{k=1}^{n}\frac{\vf_{k}}{\vt\vr Y_{k}}{\Grad p_{k}}
    =-\sum_{k=1}^{n}\frac{\vf_{k}}{\vt\vr Y_{k}}\left({\Grad p_{k}}-Y_{k}{\Grad \pi_m}\right) =\sumkN\frac{\pi_m\vf_{k}^{2}}{\vC_0\vt\vr_{k}}\geq0.
\end{equation*}

Thus, it remains to control the positive part of $\vr s(T)$ and the negative part of $(\vr s)^{0}$. From definition of the entropy (\ref{chF:entropycv})  we get
\begin{equation}\label{ch4:entform}
\vr s=\sumkN\vr Y_{k}s_{k}^{st}+\sumkN c_{v}\vr_{k}\log{{\vt}}-\sumkN \frac{\vr_{k}}{m_{k}}\log{\frac{\vr_{k}}{m_{k}}},
\end{equation}
therefore
\begin{equation*}\label{ch4:posentr}
\intO{[\vr s(T)]_{+}}\leq c\intO{\vr(T)}+c\intO{\vr\vt(T)} -\sumkN\intO{ \frac{\vr_{k}}{m_{k}}\log{\frac{\vr_{k}}{m_{k}}}(T)}.
\end{equation*}
The two first terms from the r.h.s.  are bounded due to \eqref{ch4:vrvt}, whereas to estimate the positive part of the last one we essentially use the assumption that $\Omega$ is a bounded domain. Thus, the positive part of $-x\log{x}$ is bounded  by a constant, and thus integrable over $\Omega$. $\Box$\\
In the rest of Section \ref{ch4:sect4} we show how to use Lemmas \ref{ch4:Lemma5} and \ref{ch4:Lemma6} in order to derive uniform estimates for the sequence of smooth solutions $\{\vr_N,\vu_N,\vt_N,\vr_{k,N}\}_{N=1}^\infty$ to system \eqref{ch4:1.1}.\\

\subsection{Estimates of the temperature.} One of the main consequences of \eqref{ch4:estentr} is that for $\kappa(\vr,\vt)$ satisfying \eqref{ch4:cotok} we have the following a priori estimates for the temperature
\begin{equation}\label{ch4:main-temp-est}
(1+\sqrt{\vr})\Grad\log\vt,\  (1+\sqrt{\vr})\Grad\vt^{s}\in L^{2}((0,T)\times\Omega),
\end{equation}
where $s\in [0,\frac{\alpha}{2}]$ and $\alpha\geq2$. To control the full norm of $\vt^s$ in $L^2(0,T;W^{1,2}(\Omega))$ we will apply the following version on the Korn-Poincar\'e inequality (see f.i. Theorem 10.17 in \cite{NS}):
\begin{Theorem}\label{ch4:TheoremKP}
Let $\Omega\subset \mathbb{R}^3$ be a bounded Lipschitz domain. Assume that $r$ is a non-negative function such that
$$0<M_0\leq\intO{r},\ \intO{r^\gamma}\leq K,\quad for\ a\ certain\ \gamma>1.$$
Then
$$\|\xi\|_{W^{1,p}(\Omega)}\leq C(p,M_0,K)\|\Grad\xi\|_{L^{p}(\Omega)}+
\intO{r|\xi|},$$
for any $\xi\in W^{1,p}(\Omega)$. 
\end{Theorem}

Recalling \eqref{chF:intener}, \eqref{ch4:coldp},  \eqref{ch4:vrvt}  and \eqref{ch4:main-temp-est} one can check that the assumptions of the above theorem are satisfied for $\xi=\vt$, $r=\vr$ and $p=2$. Therefore, the Sobolev imbedding gives the estimate of the norm of $\vt$ in $L^{2}(0,T; L^6(\Omega))$,
and so, due to the boundedness of $\Grad\vt^{\frac{\alpha}{2}}$ in $L^{2}((0,T) \times\Omega)$, one gets 
\begin{equation}\label{ch4:maxtemp}
\vt^{\frac{\alpha}{2}}\in L^{2}(0,T; W^{1,2}(\Omega)).
\end{equation}

\subsection{Estimates following from the Bresch-Desjardin equality}
The aim of this subsection is to derive estimates following from \eqref{ch4:first} and \eqref{ch4:both}. Summing these two expressions we obtain
\begin{multline}\label{ch4:sum}
\Dt \intO{{1\over2}\vr|\vu|^2+{1\over2}\vr|\vu+\Grad\phi(\vr)|^2}++\intO{\vS:\Grad\vu}+{1\over2}\intO{\mu(\vr)|\Grad\vu-\Grad^{T}\vu|^2}=\\
-\intO{\Grad\phi(\vr)\cdot\Grad \pi(\vt,\vr,Y)}+2\intO{\pi(\vr,\vt,Y)\Div\vu},
\end{multline}
We first need to justify that the terms  from the r.h.s. are bounded or have a negative sign so that they can be moved to the l.h.s. The main problem is to control the contribution from the molecular pressure. To this purpose we will employ estimate \eqref{ch4:estentr} coupled with the properties of the diffusion matrix $C$.\\
Denoting
\begin{equation}\label{ch4:denot}
\vC\Grad_{x_{i}}{p}=(\Grad_{x_{i}}{p})^{I},
\end{equation}
where 
\begin{equation}
{p}=\left(\begin{array}{c}
p_1\\
\vdots\\
p_n
\end{array}\right)\quad\mbox{and} \quad
\Grad{p}=\left(\begin{array}{c}
\Grad p_1\\
\vdots\\
\Grad p_n
\end{array}\right),
\end{equation}
we obtain, for every $k$-th coordinate $k\in\{1,\ldots,n\}$ and every $i$-th space coordinate $i\in\{1,2,3\}$, the following decomposition
    \begin{equation}\label{ch4:uu}
        (\Grad_{x_{i}}  p)_{k}=(\Grad_{x_{i}} {p})^{I}_{k}+\alpha_{i} Y_k.
    \end{equation}
Next, multiplying the above expression by $m_{k}$ and summing over $k\in\{1,\ldots,n\}$ one gets
    $$\alpha_{i}=\frac{\Grad_{x_{i}}(\vr\vt)}{\sum_{k=1}^{n}m_{k} Y_k}-\frac{\sum_{k=1}^{n}m_{k}(\Grad_{x_{i}}{p})^{I}_{k}}{\sum_{k=1}^{n}m_{k} Y_k }.$$
Returning \eqref{ch4:uu} we can express the full gradients of partial pressures in terms of gradients of temperature, density and the gradient of "known" part of the pressure
    \begin{equation}\label{ch4:decomp}
        \Grad{p} = ( \Grad{p})^{I} + \left(\frac{\Grad(\vr\vt)}{\sum_{k=1}^{n}m_{k} 
Y_k}-\frac{\sum_{k=1}^{n}m_{k}(\Grad{p})^{I}_{k}}{\sum_{k=1}^{n}m_{k} Y_k }\right)  Y.
    \end{equation}
As was announced, we will use the above expression in order to control the molecular part of the pressure from the r.h.s. of  \eqref{ch4:sum}.\\

\noindent {\bf {Estimate of $\Grad \pi(\vr,\vt,Y)\cdot\Grad\phi$.} }
Using definition of $\phi$ (see Lemma \ref{ch4:Lemma5}) and \eqref{chF:decompp} we obtain
    \begin{equation}\label{ch4:rs}
        \Grad\phi(\vr)\cdot\Grad \pi(\vr,\vt,Y)=\mu'(\vr)\pi'_{c}(\vr)\frac{|\Grad\vr|^2}{\vr}+
        \frac{\Grad\mu(\vr)\cdot\Grad\pi_{m}}{\vr}.
    \end{equation}
The first  term is non-negative due to \eqref{ch4:coldp}, so it can be considered on the l.h.s. of \eqref{ch4:sum} and we only need to estimate the second one. Since 
$$\Grad\pi_{m}=\sumkN(\Grad p)_{k}\quad\mbox{and} \quad\sumkN ({Y})_k=1,$$ 
we may use \eqref{ch4:decomp} to write
    \begin{multline}\label{ch4:pepi0}
        \intO{\frac{\Grad\mu(\vr)\cdot\Grad\pi_{m}}{\vr}}
        =\intO{\frac{\Grad\mu(\vr)\cdot\sumkN(\Grad p)^{I}_{k}}{\vr}}+\intO{\frac{\Grad\mu(\vr)\cdot\Grad\vr \vt}{\sumkN\vr_{k}m_{k}}}\\
        +\intO{\frac{\Grad\mu(\vr)\cdot\Grad\vt \vr}{\sumkN\vr_{k}m_{k}}}-
        \intO{\frac{\Grad\mu(\vr)\cdot\sumkN m_{k}(\Grad p)^{I}_{k}}{\sumkN\vr_{k}m_{k}}}=\sum_{i=1}^4 I_i.
    \end{multline}
Note that $I_2$ is non-negative, so we can put it to the l.h.s. of \eqref{ch4:sum}.\\
Next, $I_1$ and $I_4$ can be estimated in a similar way,  we have

\begin{equation}\label{ch4:ltop}
\intO{\frac{|\Grad\mu(\vr)||\sumkN( \Grad{p})_k^{I}|}{\vr}} \leq \ep\intO{\frac{|\Grad\mu(\vr) |^{2}\vt}{\vr}}+ c(\ep)\intO{\frac{|\sumkN(\Grad  p)^{I}_k|^{2}}{\vt\vr}},
\end{equation}
so for $\ep$ sufficiently small, the first term can be controlled  by $I_2$ thanks to  \eqref{ch4:10}.
Concerning the second integral, from (\ref{ch4:estentr}) we have

%

%
\begin{equation}\label{ch4:reduced}
\intTO{\sumkN\frac{\pi_m\vf_{k}^{2}}{\vC_0\vt\vr_{k}}}\leq c.
\end{equation}
Using \eqref{ch4:difp}, the integral may be transformed as follows
\begin{equation}\label{ch4:reduced2}
\intTO{\sumkN\frac{\vC_0(\vC\Grad  p)_k^{2}}{\pi_m\vt\vr_{k}}}\leq c,
\end{equation}
thus, due to \eqref{ch4:denot} and \eqref{ch4:assC}, the integral over time of the r.h.s. of \eqref{ch4:ltop} is bounded.
For  $I_3$ we verify that
    $$\left|\Grad\mu(\vr)\cdot\Grad\vt\frac{\vr}{\sumkN\vr_{k}m_{k}}\right|\leq c(\ep)\kappa(\vr,\vt)\frac{|\Grad\vt|^{2}}{\vt^{2}}
    +\ep\frac{\vr\vt^{2}}{\kappa(\vr,\vt)}\frac{|\Grad\mu(\vr)|^{2}}{\vr},$$
and the first term is bounded in view of \eqref{ch4:main-temp-est} whereas boundedness of the second one follows from the Gronwall inequality applied to (\ref{ch4:sum}). Indeed, note that,  due to \eqref{ch4:cotok},  $\frac{\vr\vt^{2}}{\kappa(\vr,\vt)}$ is bounded by some positive constant. \\

\noindent {\bf {Estimate of $\pi(\vr,\vt,Y)\Div\vu$.}}
By virtue of \eqref{chF:decompp} and \eqref{chF:intener} and the continuity equation
    $$\intO{\pi(\vr,\vt,Y)\Div\vu}=-\Dt\intO{\vr e_c(\vr)}+\intO{\vr\vt\left(\sumkN{\frac{Y_k}{m_k}}\right)\Div\vu}.$$
Furthermore, by the Cauchy inequality
    \begin{equation*}
        \left|\intO{\vr\vt\left(\sumkN{\frac{Y_k}{m_k}}\right)\Div\vu}\right|
        \leq c\|Y_{k}\|_{L^{\infty}(\Omega)}\left(\ep\|\sqrt{\mu(\vr)}\Div\vu\|^2_{L^{2}(\Omega)}+c(\ep)\left\|\frac{\vr\vt}{\sqrt{\mu(\vr)}}\right\|^2_{L^{2}(\Omega)}\right).
    \end{equation*}
Since $\mu(\vr)\geq \Un{\mu'}\vr$, we may write
    \begin{equation}\label{ch4:lasttop}
        \left\|\frac{\vr\vt}{\sqrt{\mu(\vr)}}\right\|_{L^{2}(\Omega)}\leq
c\|\vr\vt^2\|_{L^1(\Omega)}^{\frac{1}{2}}\leq c\|\vr \|^{\frac{1}{2}}_{L^{\frac{3}{2}}(\Omega)}\|\vt\|_{L^{6}(\Omega)}.
    \end{equation}
On account of \eqref{ch4:maxtemp},  $\vt\in L^2(0,T; L^6(\Omega))$. Moreover, the Sobolev imbedding theorem implies that 
$\|{\vr}\|_{L^{\frac{p}{2}}(\Omega)}\leq c\left\|\frac{\Grad\mu(\vr)}{\sqrt{\vr}}\right\|_{L^{2}(\Omega)}$ for $1\leq p\leq 6$,
hence the Gronwall inequality applied to \eqref{ch4:sum} implies boundedness of \eqref{ch4:lasttop}, whence the term $\ep\|\sqrt{\mu(\vr)}\Div\vu\|^2_{L^{2}(\Omega)}$ is then absorbed by $\intO{\vS:\Grad\vu}$ from the l.h.s. of \eqref{ch4:sum}.\\
Resuming, we have proved the following inequality:
\begin{multline}\label{ch4:mainest}
\mbox{ess}\sup_{t\in(0,T)}\intO{\frac{1}{2}\vr|\vu|^2+\vr e_c(\vr){1\over2}\vr|\vu+\Grad\phi(\vr)|^2(t)}+\intTO{\mu'(\vr)\pi'_{c}(\vr)\frac{|\Grad\vr|^2}{\vr}}\\
+(1-\ep)\intTO{\frac{\Grad\mu(\vr)\cdot\Grad\vr \vt}{\sumkN\vr_{k}m_{k}}}+\intTOB{\vS:\Grad\vu+{1\over2}\mu(\vr)|\Grad\vu-\Grad^{T}\vu|^2}\leq c.
\end{multline}
\noindent {\bf {Uniform estimates.}}
Taking into account all the above considerations, we can complement the so-far obtained estimates as follows\\
\begin{equation}\label{ch4:negdens}
\left\|\sqrt{{\vt}{\vr}^{-1}}\Grad\vr\right\|_{L^{2}((0,T)\times\Omega)}+\left\|\sqrt{\pi'_c(\vr)\vr^{-1}}\Grad\vr\right\|_{L^{2}((0,T)\times\Omega)}\leq c,
\end{equation}
moreover
\begin{equation}\label{ch4:posdens}
\left\|\frac{\Grad\mu(\vr)}{\sqrt{\vr}}\right\|_{L^{\infty}(0,T;L^2(\Omega))}\leq c.
\end{equation}
Concerning the velocity vector field, in addition to \eqref{ch4:vrvt} we have
    \begin{equation}
\|\sqrt{\mu(\vr)}\Grad\vu\|_{L^{2}((0,T)\times\Omega)}+\left\|{\sqrt{\mu(\vr)\vt^{-1}}}\Grad\vu\right\|_{L^{2}((0,T)\times\Omega)}\leq c.
    \end{equation}
    
 \subsection{Estimates of species densities}
Finally, we can take advantage of the entropy estimate \eqref{ch4:estentr} which together with \eqref{ch4:posdens} may be used to deduce boudedness of gradients of all species densities. 
\begin{Lemma}
For any smooth solution of \eqref{ch4:1.1} we have
\begin{equation}\label{ch4:gradcomp}
\left\|\sqrt{{1+\vt}}\Grad
\sqrt{\vr_{k}}\right\|_{L^{2}((0,T)\times\Omega)}\leq c.
\end{equation}
\end{Lemma}
\pf
First, using Remark \ref{ch4:remF} we may write
$${\frac{\pi_m\vf_k^2}{\vC_0\vt\vr_k}}=\frac{\vC_0|\Grad p_k|^2}{\pi_m\vr_k\vt}-2\frac{Y_k\vC_0\Grad p_k\cdot\Grad\pi_m}{\pi_m\vr_k\vt}+\frac{Y_k^2\vC_0|\Grad \pi_m|^2}{\pi_m\vr_k\vt}$$
which is bounded in $L^1((0,T)\times\Omega)$ on account of \eqref{ch4:reduced}. Clearly,
\begin{equation}\label{ch4:almostY}
\intTO{\frac{\vC_0|\Grad p_k|^2}{\pi_m\vr_k\vt}}\leq c\left(1+\intTO{\frac{Y_k^2\vC_0|\Grad \pi_m|^2}{\pi_m\vr_k\vt}}\right).
\end{equation}
The r.h.s. of above can be, due to \eqref{ch4:decomp}, estimated as follows
\begin{multline*}
\intTO{\frac{Y_k^2\vC_0|\Grad \pi_m|^2}{\pi_m\vr_k\vt}}=\intTO{\frac{Y_k\vC_0|\sumkN(\Grad  p)_k|^2}{\pi_m\vr\vt}}\\
\leq c\intTO{\frac{\vC_0}{\pi_m\vr\vt}\left({\left|\sumkN(\vC\Grad  p)_k\right|^2}+\frac{|\Grad(\vr\vt)|^2}{\left(\sum_{k=1}^{n}m_{k} Y_k\right)^2}+\frac{|\sum_{k=1}^{n}m_{k}(\vC\Grad{p})_{k}|^2}{\left(\sum_{k=1}^{n}m_{k} Y_k \right)^2}\right)},
\end{multline*}
which is bounded thanks to \eqref{ch4:main-temp-est}, \eqref{ch4:reduced2} and \eqref{ch4:negdens}. 
In consequence, \eqref{ch4:almostY} is bounded.
Recalling assumptions imposed on $C_0$ \eqref{ch4:assC} and the form of molecular pressure $\pi_m$, we deduce that
$$\intTO{\frac{\Un{\vC_0}(1+\vt)|\Grad\vr_k|^2}{\vr_k}}\leq c\left(1+\intTO{\frac{(1+\vt)\vr_k|\Grad\vt|^2}{\vt^2}}\right)$$
and the r.h.s. is bounded, again by  \eqref{ch4:maxY} and \eqref{ch4:main-temp-est}. $\Box$


\subsection{Additional estimates.}
In this subsection we present several additional estimates based on imbeddings of Sobolev spaces and the simple interpolation inequalities.\\

\vspace{0.2cm}

\vspace{0.2cm}

\noindent{\bf Further estimates of $\vr$.} From \eqref{ch4:coldp} and \eqref{ch4:negdens}  we deduce that there exist functions
$\xi_1(\vr)=\vr$ for $\vr<(1-\delta)$, $\xi_1(\vr)=0$ for $\vr>1$  and
$\xi_2(\vr)=0$ for $\vr<1$, $\xi_2(\vr)=\vr$ for $\vr>(1+\delta)$, $\delta>0$,
such that
$\|\Grad\xi_1^{-\frac{\gamma^-}{2}}\|_{L^{2}((0,T)\times\Omega)}$, $\|\Grad\xi_2^{\frac{\gamma^+}{2}}\|_{L^{2}((0,T)\times\Omega)}\leq c$.
Additionally in accordance to \eqref{ch4:vrvt} we are allowed to use the Sobolev imbeddings, thus
\begin{equation}\label{ch4:xipos}
\|\xi_1^{-\frac{\gamma^-}{2}}\|_{L^{2}(0,T;L^6(\Omega))},\ \|\xi_2^{\frac{\gamma^+}{2}}\|_{L^{2}(0,T;L^6(\Omega))}\leq c.
\end{equation}
\begin{Remark}
Note in particular that the first of these estimate implies that
\begin{equation}\label{ch4:!!!}
\vr(t,x)>0\quad \mbox{a.e.\ on\ }(0,T)\times\Omega.
\end{equation}
\end{Remark}
\noindent Similarly, combination of \eqref{ch4:posdens} with \eqref{ch4:vrvt} leads to
$\|\vr^{\frac{1}{2}}\|_{L^{6}(\Omega)}\leq c\left\|\frac{\Grad\mu(\vr)}{\sqrt{\vr}}\right\|_{L^{2}(\Omega)},$  
and therefore
\begin{equation}\label{ch4:maxdens}
\vr\in L^\infty(0,T;L^3(\Omega)).
\end{equation}

\vspace{0.2cm}

\noindent{\bf Estimate of the velocity vector field.} We use the
H${\rm\ddot{o}}$lder inequality  to write
\begin{equation}\label{ch4:gradu}
\|\Grad\vu\|_{L^{p}(0,T;L^{q}(\Omega))}\\
\leq c\left(1+\|\xi_1(\vr)^{-1/2}\|_{L^{2\gamma^-}(0,T;L^{6\gamma^-}(\Omega))}\right)\|\sqrt{\vr}\Grad\vu\|_{{L^{2}((0,T)\times\Omega)}},
\end{equation}
where $
p=\frac{2\gamma^-}{\gamma^-+1},\quad q=\frac{6\gamma^-}{3\gamma^-+1}.$
Therefore, Theorem \ref{ch4:TheoremKP} together with the Sobolev imbedding imply
\begin{equation}\label{ch4:estu}
\vu\in L^{\frac{2\gamma^-}{\gamma^-+1}}(0,T;L^{\frac{6\gamma^-}{\gamma^-+1}}(\Omega)).
\end{equation}
Next, by the similar argument 
\begin{equation}\label{ch4:modu}
\|\vu\|_{L^{p'}(0,T;L^{q'}(\Omega))}
\leq c\left(1+\|\xi_1(\vr)^{-1/2}\|_{L^{2\gamma^-}(0,T;L^{6\gamma^-}(\Omega))}\right)\|\sqrt{\vr}\vu\|_{{L^{\infty}(0,T;L^2(\Omega))}},
\end{equation}
with $p'=2\gamma^-,\quad q'=\frac{6\gamma^-}{3\gamma^-+1}$. By a simple interpolation between \eqref{ch4:estu} and \eqref{ch4:modu}, we obtain 
\begin{equation}\label{ch4:optu}
\vu\in L^{\frac{10\gamma^-}{3\gamma^-+3}}(0,T;L^{\frac{10\gamma^-}{3\gamma^-+3}}(\Omega)),
\end{equation}
and since $\gamma^->1$, we see in particular that $\vu\in L^{\frac{5}{3}}(0,T;L^{\frac{5}{3}}(\Omega))$.



\vspace{0.2cm}

\noindent{\bf Strict positivity of the absolute temperature.}
We now give the proof of uniform with respect to $N$ positivity of $\vt_N$.
\begin{Lemma}
Let $\{\vt_N\}_{N=1}^\infty$ be the sequence of smooth functions satisfying estimates \eqref{ch4:vrvt} and \eqref{ch4:main-temp-est}, then
\begin{equation}\label{ch4:???}
\vt_N(t,x)> 0\quad \mbox{a.e.\ on\ }(0,T)\times\Omega.
\end{equation}
\end{Lemma}
\pf The above statement is a consequence of the following estimate 
\begin{equation}\label{ch4:tempbound}
\intTO{|\log\vt_N|^{2}+|\Grad\log\vt_N|^{2}}\leq c,
\end{equation}
which can be obtained, again by application of Theorem \ref{ch4:TheoremKP} with $\xi=\log\vt_N$ and $r=\vr_N$. It remains to check that we control the $L^1(\Omega)$ norm of $\vr|\log\vt|$.
By virtue of \eqref{ch4:weakentr} we have
$$\intO{(\vr_N s_N)^{0}}\leq \intO{\vr_N s_N (T)},$$
thus substituting the form of $\vr s$ from \eqref{ch4:entform} we obtain
$$- c_v \intO{\vr_N\log{{\vt_N}}(T)}\leq\sumkN\intO{\vr_{k,N}s_{k}^{st}(T)}-\sumkN\intO{ \frac{\vr_{k,N}}{m_{k}}\log{\frac{\vr_{k,N}}{m_{k}}}(T)}-\intO{(\vr_N s_N)^{0}}$$

and the r.h.s. is bounded on account of  \eqref{ch4:maxY}, \eqref{ch4:maxdens} and the initial condition. On the other hand, the positive part of the integrant $\vr_N\log\vt_N$ is bounded from above by $\vr_N\vt_N$ which belongs to $L^\infty(0,T;L^1(\Omega))$ due to \eqref{ch4:vrvt}, so we  end up with
\begin{equation}\label{ch4:sinf}
\rm{ess}\sup_{t\in(0,T)}\intO{\left|\vr_N\log\vt_N(t)\right|}\leq c,
\end{equation}
which was the missing information in order to apply Theorem \ref{ch4:TheoremKP}.  This completes the proof of \eqref{ch4:tempbound}. $\Box$

\section{Passage to the limit}\label{ch4:sect6}
In this section we justify that it is possible to perform the limit passage in the weak formulation of system \eqref{ch4:1.1}. 
We remark that we focus only on the new features of the system, i.e. the molecular pressure and multicomponent diffusion, leaving the rest of limit passages to be performed analogously as in \cite{EZ2}. \\

\subsection { Strong convergence of the density and passage to the limit in the continuity equation.}
The strong convergence of a sequence $\{\vr_N\}_{N=1}^{\infty}$ is guaranteed by the following lemma
\begin{Lemma}\label{ch4:Lemma:vr}
If $\mu(\vr)$ satisfies \eqref{ch4:10}, then for a subsequence we have
\begin{equation}\label{ch4:pointvr}
\sqrt{\vr_N}\rightarrow\sqrt{\vr}\quad a.e.\ and\  strongly\ in\ L^2((0,T)\times\Omega).
\end{equation}
Moreover $\vr_N\rightarrow\vr$ strongly in $C([0,T];L^p(\Omega))$, $p<3$.
\end{Lemma}
For the proof see \cite{EZ2} Lemma 7.\\
In addition, from \eqref{ch4:optu}, there exists a subsequence such that
$$\vu_N\to\vu\quad\mbox{weakly\ in\ } L^{\frac{5}{3}}((0,T)\times\Omega).$$
On account of that, it is easy to let $N\to\infty$ in the continuity equation to obtain \eqref{ch4:weakcont}.

\vspace{0.2cm}
\subsection{Strong convergence of the species densities.}
Analogously we show the strong convergence of species densities. We have
%
\begin{Lemma}\label{ch4:Lemma_vrk}
Up to a subsequence the partial densities $\vr_{k,N}$, $k=1,\ldots,n$ converge strongly in $L^p(0,T;L^{q}(\Omega))$, $1\leq p<\infty, 1\leq q<3$ to $\vr_k$. In particular
\begin{equation}\label{ch4:pointvrk}
\vr_{k,N}\to\vr_k\quad a.e.\ in\  (0,T)\times \Omega.
\end{equation}
Moreover $\vr_{k,N}\to\vr_k$  in $C([0,T];L^3_{\rm{weak}}(\Omega))$.
\end{Lemma}
\pf
The estimate \eqref{ch4:gradcomp} together with \eqref{ch4:maxY} and \eqref{ch4:maxdens} give the bound for the space gradients of $\vr_{k,N}$, $k=1,\ldots,n$
\begin{equation}\label{ch4:gradvrk}
\Grad\vr_{k,N}=2\Grad\sqrt{\vr_{k,N}}\sqrt{\vr_{k,N}}\quad {\rm is\ bounded\ in\ } L^{2}(0,T;L^{3\over 2}(\Omega)).
\end{equation}
Moreover, directly from the equation of species mass conservation we obtain
\begin{equation}\label{ch4:part_dt}
\ptb{\vr_{k,N}}:=-\Div (\vr_{k,N} \vu_N)- \Div (\vf_{k,N})  +  \vr_N\vw_{k,N}
\in L^{\frac{2\alpha}{\alpha+1}}(0,T;W^{-1,\frac{6\alpha}{4\alpha+1}}(\Omega)).
\end{equation}
Indeed, the most restrictive term is the diffusion flux, which can be rewritten as
\begin{equation}\label{ch4:fineF}
{\vf_{i,N}}=
-\frac{\vC_0}{\pi_{m,N}}\left(\Grad\vt_N\frac{\vr_{i,N}}{m_i}+\Grad\vr_{i,N}\frac{\vt_N }{m_i}
-{Y_{i,N}}\sumkN \Grad\vt_N\frac{\vr_{k,N}}{m_k}-{Y_{i,N}}\sumkN\Grad\vr_{k,N}\frac{\vt_N}{m_k}
\right).
\end{equation}
Due to \eqref{ch4:maxY} we have that
$$\frac{\vC_0}{\pi_{m,N}}\left|\Grad\vt_N\frac{\vr_{i,N}}{m_i}-{Y_{i,N}}\sumkN \Grad\vt_N\frac{\vr_{k,N}}{m_k}\right|\leq c\left(|\Grad\log\vt_N|+|\Grad\vt_N|\right)\vr_N,$$
which is bounded in $L^2(0,T;L^{\frac{3}{2}})$ on account of\eqref{ch4:main-temp-est} and \eqref{ch4:maxdens}. Similarly,
\begin{equation*}
\frac{\vC_0}{\pi_{m,N}}\left|\Grad \vr_{i,N}\frac{\vt_N}{m_i}-{Y_{i,N}}\sumkN\Grad \vr_{k,N} \frac{\vt_N}{m_k}
\right|
\leq c \sqrt{(\vt_N+1)\vr_N}\sumkN\left| \sqrt{\vt_N+1}\Grad \sqrt{\vr_{k,N}}\right|,
\end{equation*}
thus, according to \eqref{ch4:gradcomp} it remains to control the norm of $\sqrt{(\vt_N+1)\vr_N}$ in $L^{p}((0,T)\times\Omega)$ for some $p>2$. By \eqref{ch4:maxtemp} and \eqref{ch4:maxdens} we deduce that $\vr_N\vt_N\in L^\alpha(0,T; L^{\frac{3\alpha}{\alpha+1}}(\Omega))$, so \eqref{ch4:part_dt} is verified.
In this manner we actually proved that the sequence of functions
$$\{t\to\intO{{\vr_{k,N}} \phi }\},\quad \phi\in C^\infty _c(\Omega)$$
is uniformly bounded and equicontinuous in $C([0,T])$, hence, the Arzel\'a-Ascoli theorem yields
$$\intO{{\vr_{k,N}} \phi }\to \intO{{\vr_{k}} \phi }\quad\mbox{in}\ C([0,T]).$$
Since $\vr_{k,N}$ is bounded in $L^\infty(0,T;L^3(\Omega))$ and due to density argument, this convergence extends to each $\phi \in L^{\frac{3}{2}}(\Omega)$.\\
Finally, the Aubin-Lions argument implies  the strong convergence of the sequence $\vr_{k,N}$ to $\vr_k$ in $L^p(0,T;L^{q}(\Omega))$ for $p=2, q<3$, but due to \eqref{ch4:maxY} and \eqref{ch4:maxdens} it can be extended to the case $p<\infty.$ $\Box$

\vspace{0.2cm}

\subsection{Strong convergence of the temperature.}
From estimate \eqref{ch4:maxtemp} we deduce existence of a subsequence such that
\begin{equation}\label{ch4:conv_teta}
\vt_N\to\vt\quad\mbox{weakly\ in\ }L^2(0,T;W^{1,2}(\Omega)),
\end{equation}
however, time-compactness cannot be proved directly from the internal energy equation \eqref{ch4:tempeq}. The reason for this is lack of control over a part of the heat flux proportional to $\vr\vt^{\alpha}\Grad\vt$. This obstacle can be overcome by deducing analogous information from the entropy equation \eqref{ch4:ent}.
\vspace{0.2cm}

We will first show that all of the terms appearing in the entropy balance \eqref{ch4:ent} are nonnegative or belong to $W^{-1,p}((0,T)\times\Omega)$, for some $p>1$. \\
Indeed, first recall that due to \eqref{ch4:entform} 
\begin{equation*}
|\vr_N s_N|\leq c\left(\vr_N+\vr_N|\log{\vt_N}|+\sumkN \vr_{k,N}|\log{\vr_{k,N}}|\right)
\end{equation*}
and
\begin{equation*}
|\vr_N s_N\vu_N|\leq c\left(|\vr_N\vu_N|+|\vr_N\log{\vt_N}\vu_N|+\sumkN |\vr_{k,N}\log{\vr_{k,N}}\vu_N|\right),
\end{equation*}
whence due to  \eqref{ch4:vrvt}, \eqref{ch4:maxdens} and \eqref{ch4:tempbound} we deduce that
\begin{equation}\label{ch4:2ent}
\{\vr_N s_N\}_{N=1}^\infty\quad \mbox{ is\ bounded\ in}\ L^{2}((0,T)\times\Omega),
\end{equation}
moreover
\begin{equation}\label{ch4:beg_ent_flux}
\{\vr_N s_N\vu_N\}_{N=1}^\infty\quad \mbox{ is\ bounded\ in}\ L^{2}(0,T; L^{\frac{6}{5}}(\Omega)).
\end{equation}

The entropy flux is due to \eqref{chF:defg} and \eqref{chF:cotoQ} equal to
$$\frac{\vQ}{\vt}-\sumkN \frac{g_{k}}{\vt}\vf_{k}=\frac{\kappa(\vr,\vt)\Grad\vt}{\vt}+\sumkN s_k\vf_{k}.$$
The first part can be estimated as follows
$$\left|\frac{\kappa(\vr_N,\vt_N)\Grad\vt_N}{\vt_N}\right|\leq |\Grad\log\vt_N|+
|\vr_N\Grad\log\vt_N|+|\vt_N^{\alpha-1}\Grad\vt_N|+|\vr_N\vt_N^{\alpha-1}\Grad\vt_N|,$$
where the most restrictive term can be controlled as follows $
|\vr_N\vt_N^{\alpha-1}\Grad\vt_N|\leq|\sqrt{\vr_N}\vt_N^{\frac{\alpha}{2}}||\sqrt{\vr_N}\Grad\vt_N^{\frac{\alpha}{2}}|$, 
which is bounded  on account of \eqref{ch4:main-temp-est} provided ${\vr_N}\vt_N^{\alpha}$ is bounded in $L^{p}((0,T)\times\Omega)$ for $p>1$, uniformly with respect to $N$. Note that for $0\leq \beta\leq 1$ we have
${\vr_N}\vt_N^{{\alpha}}=(\vr_N\vt_N)^{\beta}\vr_N^{1-\beta}\vt_N^{\alpha-\beta}$,
where $(\vr_N\vt_N)^{\beta}$, $\vr_N^{1-\beta}$, $\vt_N^{\alpha-\beta}$ are uniformly bounded in 
$L^{\infty}(0,T;L^{\frac{1}{\beta}}(\Omega))$,  $L^{\infty}(0,T;L^{\frac{3}{1-\beta}}(\Omega))$ and $L^{\frac{\alpha}{\alpha-\beta}}(0,T;L^{\frac{3\alpha}{\alpha-\beta}}(\Omega))$, respectively.
Therefore
\begin{equation}\label{ch4:hardheat}
\left\{\frac{\kappa(\vr_N,\vt_N)\Grad\vt_N}{\vt_N}\right\}_{N=1}^\infty\quad \mbox{ is\ bounded\ in}\ L^{p}(0,T; L^{q}(\Omega)),
\end{equation}
for $p$ and $q$ satisfying
$\frac{1}{p}=\frac{\alpha-\beta}{\alpha},\quad \frac{1}{q}=\beta+\frac{1-\beta}{3}+\frac{\alpha-\beta}{3\alpha}.$
In particular $p,q>1$ provided
$0<\beta<\frac{\alpha}{2\alpha-1}.$\\

The remaining part of the entropy flux is equal to
\begin{equation*}\label{ch4:entdiff}
\sumkN s_{k,N}\vf_{k,N}=\sumkN\frac{\vf_{k,N}}{m_k}+c_v\sumkN{\log\vt_N\vf_{k,N}}
-\sumkN\frac{\vf_{k,N}}{m_k}\log\frac{\vr_{k,N}}{m_k},
\end{equation*}
where the middle term vanishes due to \eqref{ch4:zerototdiff}. The worst term to estimate is thus the last one, we rewrite it using \eqref{ch4:difp} in the following way
\begin{equation*}
-\frac{\vf_{i,N}}{m_i}\log\frac{\vr_{i,N}}{m_i}
=\frac{\vC_0\Grad p_{i,N}}{\pi_{m,N} m_i}\log\frac{\vr_{i,N}}{m_i}-\frac{\vr_{i,N}}{\vr_N}\frac{\vC_0\Grad\pi_{m,N}}{\pi_{m,N} m_i}\log\frac{\vr_{i,N}}{m_i}
\end{equation*}
for $i=1,\ldots,n$. Both parts have the same structure, so we focus only on the first one, we have
\begin{multline*}
\frac{\vC_0}{\pi_{m,N}}\left|\frac{\Grad p_{i,N}}{m_i}\log\frac{\vr_{i,N}}{m_i}\right|
\leq c|\sqrt{(\vt_N+1)\vr_{i,N}}\log\vr_{i,N}||\sqrt{\vt_N+1}\Grad\sqrt{\vr_{i,N}}|\\
+c\sqrt{\vr_N}\left(|\Grad\log\vt_N|+|\Grad\vt_N|\right)|\sqrt{\vr_{i,N}}\log\vr_{i,N}|.
\end{multline*}
Using \eqref{ch4:main-temp-est},  \eqref{ch4:maxtemp}, \eqref{ch4:gradcomp} and \eqref{ch4:maxdens} we finally arrive at
\begin{equation}\label{ch4:hardentr}
\left\{\sumkN s_{k,N}\vf_{k,N}\right\}_{N=1}^\infty\quad \mbox{ is\ bounded\ in}\ L^{p}((0,T)\times(\Omega)),\  \mbox{for}\  1<p<\frac{4}{3}.
\end{equation}

%
%

We are now ready to proceed with the proof of strong convergence of the temperature. To this end we will need the following  variant of the Aubin-Lions Lemma.
\begin{Lemma}\label{ch4:LLions}
Let $g^{N},\ h^{N}$ converge weakly to $g,\ h$ respectively in $L^{p_{1}}(0,T;L^{p_{2}}(\Omega))$,  $L^{q_{1}}(0,T;L^{q_{2}}(\Omega))$, where $1\leq p_{1},\ p_{2}\leq\infty$ and
\begin{equation}\label{ch4:cond0}
\frac{1}{p_{1}}+\frac{1}{q_{1}}=\frac{1}{p_{2}}+\frac{1}{q_{2}}=1.
\end{equation}
Let us assume in addition that
\begin{equation}\label{ch4:condI}
\frac{\partial g^{N}}{\partial t}\ is\ bounded\ in\ L^{1}(0,T;W^{-m,1}(\Omega))\ for\ some\ m\geq0\ independent\ of\ N
\end{equation}
\begin{equation}\label{ch4:condII}
\|h^{N}-h^{N}(\cdot+\xi,t)\|_{L_{q_{1}}(L_{q_{2}})}\rightarrow 0\ as\ |\xi|\rightarrow 0,\ uniformly\ in\ N.
\end{equation}
Then $g^{N}h^{N}$ converges to $gh$ in the sense of distributions on $\Omega\times(0,T)$.
\end{Lemma}
For the proof see \cite{PLL}, Lemma 5.1.\\

Taking $g^N=\vr_N s_N$ and $h^N=\vt_N$ we verify, due to \eqref{ch4:2ent}, \eqref{ch4:maxtemp} and \eqref{ch4:conv_teta}
, that 
 conditions \eqref{ch4:cond0}, \eqref{ch4:condII} are satisfied with $p_1,p_2,q_1,q_2=2.$  Moreover, for $m$ sufficiently large $L^1(\Omega)$ is  imbedded into $W^{-m,1}(\Omega)$, thus by the previous considerations, condition \eqref{ch4:condI} is also fulfilled. Therefore, passing to the subsequences we may deduce that $\lim_{N\to\infty}\vr_N s(\vr_N,\vt_N,Y_N)\vt_N=\Ov{\vr s(\vr,\vt,Y)}\vt.$ 
On the other hand, $\vr_N$ converges  to $\vr$ a.e. on $(0,T)\times\Omega$, hence $\Ov{\vr s(\vr,\vt,Y)}\vt=\vr\Ov{ s(\vr,\vt,Y)}\vt,$
in particular, we have that
\begin{multline}\label{ch4:equalt}
\sumkN\Ov{\frac{\vr_k}{m_k}\vt}+c_v\Ov{\vr\log\vt\ \vt}-\sumkN\Ov{\frac{\vr_k}{m_k}\log \frac{\vr_k}{m_k}\vt}=
\sumkN\Ov{\frac{\vr_k}{m_k}}\vt+c_v\vr\Ov{\log\vt}\ \vt-\sumkN\Ov{\frac{\vr_k}{m_k}\log \frac{\vr_k}{m_k}}\vt.
\end{multline}
Combining  Lemma \ref{ch4:Lemma_vrk} with \eqref{ch4:conv_teta} we identify 
\begin{equation*}
\sumkN\Ov{\frac{\vr_k}{m_k}\vt}-\sumkN\Ov{\frac{\vr_k}{m_k}\log \frac{\vr_k}{m_k}\vt}=\sumkN{\frac{\vr_k}{m_k}\vt}-\sumkN{\frac{\vr_k}{m_k}\log \frac{\vr_k}{m_k}\vt},
\end{equation*}
so  \eqref{ch4:equalt} implies that 
$\vr\Ov{\log\vt\vt}=\vr\Ov{\log\vt}\vt.$
This in turn yields that
$\Ov{\log\vt\vt}= \log\vt\vt\quad$ {a.e.\ on} $(0,T)\times\Omega$,  
 since $\vr> 0$ a.e. on $(0,T)\times\Omega$, which, due to convexity of function $x\log x$, gives rise to
 \begin{equation}\label{ch4:pointvt}
\vt_N\to\vt\quad \mbox{a.e.\ on} \ (0,T)\times\Omega.
\end{equation}

\subsection{Limit in the momentum equation, the species mass balance equations and the global total energy balance}
Having proven pointwise convergence of sequences $\{\vr_N\}_{N=1}^\infty$, $\{\vr_{k,N}\}_{N=1}^\infty$ and $\{\vt_N\}_{N=1}^\infty$ we are ready to perform the limit passage in all the nonlinear terms appearing in the momentum equation, the species mass balance equations and the total global energy balance.\\

\vspace{0.2cm}

\noindent{\bf (i) Limit in the convective term.} Estimate \eqref{ch4:estu} implies that for $0\leq\epsilon\leq1/2$ we have
\begin{equation}\label{ch4:convective}
\|\sqrt{\vr}\vu\|_{L^{p'}(0,T;L^{q'}(\Omega))}
\leq\|\vr\|^{1/2-\epsilon}_{ L^{\infty}(0,T;L^{3}(\Omega))}\|\sqrt{\vr}\vu\|_{ L^{\infty}(0,T;L^{2}(\Omega))}^{2\epsilon}\|\vu\|_{L^{\frac{2\gamma^-}{\gamma^-+1}}(0,T;L^{\frac{6\gamma^-}{\gamma^-+1}}(\Omega))}^{1-2\epsilon},
\end{equation}
where $p',\ q'$ are given by
$\frac{1}{p'}=\frac{1-2\epsilon}{\frac{2\gamma^-}{\gamma^-+1}},\quad \frac{1}{q'}=\frac{1/2-\epsilon}{3}+\frac{2\epsilon}{2}+\frac{1-2\epsilon}{\frac{6\gamma^-}{\gamma^-+1}}.$ 
Taking $\epsilon>\frac{1}{2(\gamma^-+1)}$ we have $p',q'>2$, provided $\gamma^->1$, so the convective term converges weakly to $\Ov{\vr\vu\otimes\vu}$ in $L^p((0,T)\times\Omega)$ for some $p>1$. To identify the limit, we prove the following lemma.
\begin{Lemma}
Let $p>1$, then up to a subsequence we have
\begin{eqnarray}
\vr_N\vu_N&\rightarrow& \vr\vu\quad   in\   C([0,T];L^\frac{3}{2}_{\rm{weak}}(\Omega))\label{ch4:Cpweak},\\
\vr_N\vu_N\otimes\vu_N&\rightarrow& \vr\vu\otimes\vu\quad  weakly\  in\   L^{p}((0,T)\times\Omega)\label{ch4:convconv}.
\end{eqnarray}
\end{Lemma}
\pf  We already know that $\vr_N$ converges  to $\vr$ a.e. on $(0,T)\times\Omega$. Moreover, due to \eqref{ch4:estu}, up to extracting a subsequence, $\vu_N$ converges weakly to $\vu$ in  $L^p(0,T;L^q(\Omega))$ for $p>1,\ q>3$. Therefore, the uniform boundedness of the sequence $\vr_N\vu_N$ in $L^{\infty}(0,T; L^{\frac{3}{2}}(\Omega))$ implies that 
$$\vr_N\vu_N \rightarrow \vr\vu\quad \mbox{ weakly}^*\mbox{\  in}\   L^{\infty}(0,T;L^{\frac{3}{2}}(\Omega)).$$
Now, we are aimed at improving the time compactness of this sequence. Using the momentum equation, we show that the sequence of functions
$$\{t\to\intO{\vr_N\vu_N\phi }\}_{N=1}^\infty$$
is uniformly bounded and equicontinuous in $C([0,T])$, where $\phi\in C^\infty_c(\Omega)$. But since the smooth functions are dense in $L^3(\Omega)$, applying the Arzel\`a-Ascoli theorem, we show \eqref{ch4:Cpweak}.

On the other hand, $\vu_N$ is uniformly bounded in $L^p(0,T;W^{1,q}(\Omega))$ for $p>1,\ q>\frac{3}{2}$, so it converges to $\vu$ weakly in this space. Since $W^{1,q}(\Omega)$, $q>\frac{3}{2}$ is compactly embedded into $L^3(\Omega)$, by \eqref{ch4:Cpweak}, we obtain \eqref{ch4:convconv}. $\Box$

\vspace{0.2cm}

\noindent{\bf (ii) Limit in the stress tensor.} 
\begin{Lemma}
If $\mu(\vr)$, $\nu(\vr)$ satisfy \eqref{ch4:10}, then for a subsequence we have
\begin{equation}\label{ch4:conv_stress}
\begin{array}{l}
\mu(\vr_N)\D(\vu_N)\rightarrow\mu(\vr)\D(\vu)\quad \mbox{weakly\ in}\ L^{p}((0,T)\times\Omega)\\
\nu(\vr_N)\Div\vu_N\rightarrow\nu(\vr)\Div\vu\quad\mbox{weakly\ in}\  L^{p}((0,T)\times\Omega)
\end{array}
\quad for\ p>1.
\end{equation}
\end{Lemma}
\pf Due to \eqref{ch4:gradu}, there exists a subsequence such that 
$$\Grad\vu_N\to\Grad\vu \quad \mbox{weakly\ in\ } L^p(0,T;L^q(\Omega))\ \mbox{for}\  p>1, q>\frac{3}{2}.$$  
Moreover $\mu(\vr_N),\nu(\vr_N)$ are bounded in $L^\infty(0,T;L^3(\Omega))$, on account of \eqref{ch4:10}. Thus, \eqref{ch4:conv_stress} follows by application of Lemma \ref{ch4:Lemma:vr}. $\Box$\\

\vspace{0.2cm}

\noindent{\bf (iii) Strong convergence of the cold pressure.} It follows from estimates
\eqref{ch4:vrvt} combined with \eqref{ch4:xipos} and the Sobolev imbedding theorem that
\begin{equation}\label{ch4:5/3}
\|\pi_c(\vr_{N})\|_{L^{5\over3}((0,T)\times\Omega)}\leq \|\pi_c(\vr_N)\|_{L^{\infty}(0,T;L^{1}(\Omega))}^{2\over5}\|\pi_c(\vr_N)\|_{L^{1}(0,T;L^{3}(\Omega))}^{3\over5}\leq c. 
\end{equation}
Having this, strong convergence of $\vr_N$  implies convergence of $\pi_c(\vr_n)$ to $\pi_c(\vr)$ strongly in $L^{p}((0,T)\times\Omega)$ for $1\leq p<\frac{5}{3}$.\\

\vspace{0.2cm}

\noindent{\bf (iv) Convergence of the diffusion terms.} In the proof of Lemma \eqref{ch4:Lemma_vrk} it was shown in particular that
$\left\{\vf_{k,N}\right\}_{N=1}^\infty$ { is bounded in} $L^{\frac{4}{3}}((0,T)\times(\Omega)).$ 
By the weak convergence of  $\Grad\vr_{k,N}$, $\Grad\vt_N$ to $\Grad\vr_k$, $\Grad\vt$, respectively, deduced  from \eqref{ch4:gradvrk} and \eqref{ch4:conv_teta} together with \eqref{ch4:pointvrk}, \eqref{ch4:pointvt} and \eqref{ch4:!!!} we check that it is possible to let $N\to \infty$ in all terms of \eqref{ch4:fineF}. In other words, we have 
$$\vf_{k}(\vr_N,\vt_N,\vr_{k,N})\to\vf_k(\vr,\vt,\vr_k)\quad\mbox{weakly\ in}\  L^{\frac{4}{3}}((0,T)\times\Omega), \ k\in\{1,\ldots,n\}.$$

The convergence results established above are sufficient to perform the limit passage in the momentum, the total global energy balance and the species mass balance equations and to validate, that the limit quantities satisfy the weak formulation \eqref{ch4:weakmom},\eqref{ch4:weaktoten} and \eqref{ch4:weakspec}.
\subsection{Limit in the entropy inequality}
In view of (\ref{ch4:2ent}-\ref{ch4:hardentr}) and the remarks from the previous subsection, it is easy to pass to the limit $N\to\infty$ in all terms appearing in \eqref{ch4:weakentr}, except the entropy production rate $\sigma$. \\
However, in accordance with \eqref{ch4:estentr} we still have that
$$\left\{\sqrt{\frac{\mu(\vr_N)}{\vt_N}}\left(\Grad\vu_N+(\Grad\vu_N)^T-\frac{2}{3}\Div\vu_N\right)\right\}_{N=1}^\infty\quad\mbox{is\ bounded\ in\ }L^2((0,T)\times\Omega).$$
Moreover, by virtue of \eqref{ch4:gradu}, \eqref{ch4:pointvr} and \eqref{ch4:pointvt} we deduce
$$\sqrt{\frac{\mu(\vr_N)}{\vt_N}}\left(\Grad\vu_N+(\Grad\vu_N)^T-\frac{2}{3}\Div\vu_N\right)\to\sqrt{\frac{\mu(\vr)}{\vt}}\left(\Grad\vu+(\Grad\vu)^T-\frac{2}{3}\Div\vu\right)$$
weakly in $L^2((0,T)\times\Omega).$ Evidently, we may treat all the remaining terms
$$ \left\{\sqrt{\frac{\frac{2}{3}\mu(\vr_N)+\nu(\vr_N)}{\vt_N}}\Div\vu_N\right\}_{N=1}^{\infty},\left\{\frac{\sqrt{\kappa(\vr_N,\vt_N)}}{\vt_N}\Grad\vt_N\right\}_{N=1}^\infty,\quad \left\{\frac{\sqrt{\pi_m(\vt_N,Y_{N})}}{\sqrt{\vC_0\vr_{k,N}\vt_N}}\vf_{k,N}\right\}_{N=1}^\infty$$
in the similar way using the fact that they are linear with respect to the weakly convergent sequences of gradients of $\vu_N$, $\vt_N$ and $\vr_{k,N}$.
Thus, preserving the sign of the entropy inequality \eqref{ch4:weakentr} in the limit $N\to\infty$ follows by the lower semicontinuity of convex superposition of operators.\\

Our ultimate goal is to show  that the limit entropy $\vr s$ attains its initial value at least in the weak sense. We have the following result
\begin{Lemma} Let $\vr,\vu,\vt,\vr_1,\ldots,\vr_n$ be a weak variational entropy solution to \eqref{ch4:1.1} in the sense of Definition \ref{ch4:def1}. Then the entropy $\vr s$ satisfies
\begin{equation}\label{ch4:wel}
\rm{ess}\lim_{\tau\to0^+}\intO{\left(\vr s\right)(\tau)\phi}\to \intO{(\vr s)^0\phi},\quad \forall\phi\in C^\infty(\Omega).
\end{equation}
\end{Lemma}
\pf As a consequence of \eqref{ch4:weakentr} we know that
$$\intO{(\vr s(\vt,\vr_k))(\tau^+)\phi}\geq\intO{(\vr s(\vt,\vr_k))(\tau^-)\phi}$$
where $\phi\in C^\infty(\Omega)$, $\phi\geq0$ and $(\vr s(\vt,\vr_k))(\tau^+)\in{\cal M}^+(\Omega)$, $\tau\in[0,T)$, $(\vr s(\vt,\vr_k))(\tau^-)\in{\cal M}^+(\Omega)$, $\tau\in(0,T]$ are the one sided limits of $\vr s(\tau)$. Note that due to \eqref{ch4:sinf} $\vr s\in L^\infty(0,T;L^1(\Omega))$, thus for any Lebesgue point of $\tau\mapsto\vr s(\tau,\cdot)$ these signed measures coincide with a function $\vr s(\tau,\cdot)\in L^1(\Omega)$ which satisfies

\begin{multline}\label{ch4:intmeas}
\intO{\vr s(\tau)\phi}-{\left\langle\sigma,\phi \right\rangle}\\=
\intO{(\vr s)^0\phi}-\int_{0}^\tau\intO{\vr s\vu\cdot\Grad\phi}~\dt+\int_{0}^\tau\intO{\left(\frac{\vQ}{\vt}-\sumkN\frac{g_k}{\vt}\vf_k\right)\cdot\Grad\phi}~\dt,
\end{multline}
for any test function $\phi\in C^\infty([0,\tau]\times\Omega)$, $\phi\geq0$, where $\sigma\in{\cal M}^+([0,T]\times\Omega)$ is a nonnegative measure such that
$$\sigma\geq
{\frac{\vS:\Grad\vu}{\vt}
-{\vQ\cdot\Grad\vt\over{\vt^{2}}}-\sumkN\frac{\vf_k}{m_k}\cdot\Grad\left(\frac{g_k}{\vt}\right)
-{\sumkN g_k \vr\vw_k}}.$$
In order to show \eqref{ch4:wel} we thus need to justify that $\sigma$ is absolutely continuous with respect to Lebesgue measure on $[0,\tau]\times\Omega$. To this end we use in \eqref{ch4:intmeas} a test function $\phi=\vt^0$, we get
\begin{equation}\label{ch4:A}
\intO{\vr s(\tau)\vt^0}-{\left\langle\sigma,\vt^0\right\rangle}=
\intO{(\vr s)^0\vt^0}+\int_{0}^\tau\intO{\vc{H}\cdot\Grad\vt^0}~\dt,
\end{equation}
where, on account of (\ref{ch4:beg_ent_flux}-\ref{ch4:hardentr})
\begin{equation}\label{ch4:propH}
\vc{H}=\vr s\vu+\frac{\vQ}{\vt}-\sumkN\frac{g_k}{\vt}\vf_k\in L^{p}((0,T)\times\Omega),\ \mbox{for\ some}\ p>1.
\end{equation}
Now,  testing  \eqref{ch4:weaktoten} with  $\phi_m\in C^\infty[0,T)$ such that $\phi_m\to 1$ pointwisely in $[0,\tau)$, $\phi_m\to 0$ pointwisely in $[\tau,T)$, $0< \tau<T$ and passing to the limit with $m$, we obtain
\begin{equation}\label{ch4:B}
\intOB{\frac{\left|(\vr\vu)^0\right|^2}{2\vr^0}+\vr^0 e(\vr^0,\vt^0,\vr_1^0,\ldots,\vr_n^0)}=\intO{\left(\frac{|\vr\vu|^2}{2\vr}+\vr e(\vr,\vt,\vr_1,\ldots,\vr_n)\right)(\tau)}.
\end{equation}
Combining \eqref{ch4:A} with \eqref{ch4:B}, we thus get
\begin{multline}\label{ch4:this}
\intOB{\frac{|\vr\vu|^2}{2\vr}(\tau)-\frac{\left|(\vr\vu)^0\right|^2}{2\vr^0}}+\intO{\underbrace{\left[\vr e(\tau)-\vt^0\vr s(\tau)\right]-\left[\vr^0 e^0-\vt^0(\vr s)^0\right]}_{I^*}}+{\left\langle\sigma,\vt^0\right\rangle}\\
=-\int_{0}^\tau\intO{\vc{H}\cdot\Grad\vt^0}~\dt.
\end{multline}
By the Fatou lemma
$$\rm{ess}\liminf_{\tau\to0^+}\intOB{\frac{|\vr\vu|^2}{2\vr}(\tau)-\frac{\left|(\vr\vu)^0\right|^2}{2\vr^0}}\geq0,$$
in addition

{$$\lim_{\tau\to0^+}\int_{0}^\tau\intO{\vc{H}\cdot\Grad\vt^0}~\dt=0,$$}
on account of \eqref{ch4:propH}. Moreover, recalling \eqref{chF:enerdef} and \eqref{chF:entropycv}, we recast $I^*$ as follows
\begin{multline}
I^*=\sumkN(e_k^{st}-\vt^0s_k^{st})(\vr_k(\tau)-\vr_k^0)+\vt^0\sumkN\left[\frac{\vr_k}{m_k}\log\frac{\vr_k}{m_k}(\tau)-\frac{\vr^0_k}{m_k}\log\frac{\vr^0_k}{m_k}\right]+\left[\vr e_c(\vr)(\tau)-\vr^0e_c(\vr^0)\right]\\
+c_v\left[\vr\vt(\tau)-\vr^0\vt^0-\vr\log\vt(\tau)+\vr^0\log\vt^0\right]=\sum_{i=1}^4 I_i.
\end{multline}
Evidently $\intOB{I_1+I_2+I_3}\to0$ in view of weak continuity of $\vr$ and $\vr_k$, $k=1,\ldots,n$. Concerning the last term, we have
$$I_4=c_v\underbrace{\vr(\tau)\left[\vt(\tau)-\log\vt(\tau)-\vt^0+\log\vt^0\right]}_{\geq0}+c_v(\vr(\tau)-\vr^0)(\vt^0-\log\vt^0),$$
therefore $\lim_{\tau\to0^+}\intO{I_4}\geq 0$. 
Since the entropy production rate is always nonnegative, \eqref{ch4:this} together with above remarks yields
$\rm{ess}\lim_{\tau\to0^+}{\left\langle\sigma,\vt^0\right\rangle}=0$,
whence
$$\rm{ess}\lim_{\tau\to0^+}\sigma[[0,\tau]\times\Omega]=0.~\Box$$

\medskip 

{\footnotesize {\bf Acknowlegdement.} The author was supported by the International Ph.D.
Projects Programme of Foundation for Polish Science operated within the
Innovative Economy Operational Programme 2007-2013 funded by European
Regional Development Fund (Ph.D. Programme: Mathematical
Methods in Natural Sciences)  and partly supported by Polish MN grant IdP2011000661.}

\end{document}